%% file: arxiv.tex
\documentclass[titlepage,11pt,final]{article}
\input mac.tex

\usepackage[pdfencoding=auto, psdextra, colorlinks,citecolor=blue!90!black, urlcolor= black]{hyperref}

\usepackage{enumitem}
\usepackage{bm}
\usepackage{mathrsfs}
\usepackage{framed}
\newcommand{\dd}{{\rm d}}

\usepackage{tcolorbox}

\newcommand{\dFM}[1]{d_{{\rm FM(}#1{\rm )}}}

\begin{document}


\begin{center}
\begin{large}
{\bf Approximating Rockafellians Mitigate Distributional Perturbations:\\ Discontinuous Integrands and Chance-Constrained Applications}
\smallskip
\end{large}
\vglue 0.7truecm
\begin{tabular}{ccc}
  \begin{large} {\sl Lai Tian} \end{large}& & \begin{large} {\sl Johannes O. Royset} \end{large}\\
  Daniel J. Epstein Department of & & Daniel J. Epstein Department of \\
  Industrial and Systems Engineering & & Industrial and Systems Engineering\\
  University of Southern California & & University of Southern California\\ 
  laitian@usc.edu & & royset@usc.edu 
\end{tabular}

\vskip 0.2truecm

\end{center}

\vskip 0.3truecm
\noindent {\bf Abstract}. \quad  In this paper, we show how approximating Rockafellians serve as a principled and effective alternative for improving the stability of stochastic programs under distributional changes. Unlike previous efforts that focus on special distributions and continuous integrands, our results accommodate general probability distributions and discontinuous integrands. Thus, our results apply to chance-constrained programs, for which we obtain improved qualitative and quantitative stability results under weaker assumptions pertaining to metric subregularity and upper outer-Minkowski content.

\vskip 0.2truecm

\halign{&\vtop{\parindent=0pt
   \hangindent2.5em\strut#\strut}\cr
                

\vskip 0.2truecm

{\bf Date}:\quad \ July 21, 2025 \cr}

\baselineskip=15pt

\section{Introduction}\label{sec:intro}

Solutions of a stochastic optimization problem tend to change disproportionally under small changes to its probability distribution. This is especially troubling because it is practically impossible to determine the ``right'' distribution in a real-world application. In this paper, we show that careful adjustments to a stochastic optimization problem can significantly improve the stability of its solutions under distributional perturbations. The adjustments amount to Rockafellian relaxations related to those in \cite{RoysetChenEckstrand.24,antil2024rockafellian,deride2024approximations}. In contrast to those earlier efforts and their focus on finite, discrete, or continuous distributions, we consider nearly arbitrary (Borel) probability distributions and address discontinuous integrands including the kind arising in chance constraints. Specifically,
for a closed set $\Xi\subseteq \reals^d$ and a Borel probability distribution $\mu$ on $\Xi$, we center on the {\em actual problem} 
\begin{equation}\label{eq:general-problem}
	\nnmin_{x\in \reals^n}\ \phi(x) = g_0(x) + h\big(\Ex_\mu[ G(\bm{\xi},x) ]\big),
\end{equation}
where $g_0:\reals^n \to \Reals := [-\infty,\infty]$ is a proper, lsc function, $h:\reals^m \to \Reals$ is a proper, lsc, nondecreasing function, and $G:\Xi \times \reals^n \to \reals^m$ is a vector-valued mapping with random lsc component functions satisfying certain boundedness properties clarified later. 
Among many other possibilities, a specific setting of \eqref{eq:general-problem} motivating our development has $h$ as the indicator function of $(-\infty,0]^m$ and $G(\xi,x)=(g_1(\xi,x),\ldots, g_m(\xi, x))$ with component function
$g_i:\Xi \times \reals^n \to \{b_i - 1, b_i\}$ and $b_i \in [0,1]$ for each $i \in \{1,\ldots, m\}$.
In this setting,  the function $x\mapsto h\big(\Ex_\mu[ G(\bm{\xi},x) ]\big)$ can be interpreted as imposing certain chance constraints (see Section~\ref{sec:cc} for details).
We are interested in designing adjustments to \eqref{eq:general-problem} that are stable to perturbations of the probability distribution $\mu$, possibly caused by an adversary.

 The systematic study of how stochastic problems respond to distributional changes can be traced back to early works in the field \cite{klatte1987note,robinson1987stability,kall1987approximations,romisch1991distribution,romisch1991stability}; see also the surveys \cite{schultz2000some,romisch2003stability} and the books \cite{Rachev91,ShapiroDentchevaRuszczynski.21}.
For chance-constrained programs, the works \cite{romisch1991stability,henrion1999metric,rachev2002quantitative,henrion2004holder,van2024inner} propose regularity conditions ensuring stability (in some sense) to distributional perturbations for various settings. 
These regularity conditions are usually related to continuity properties of the feasible set mapping, for example ensured by metric regularity, but these could be limiting for certain applications.

One general approach to tackle distributional uncertainty is to adopt a conservative viewpoint, which gives rise to risk measures and distributionally robust formulations; see, e.g., \cite{royset2025risk,chen2010cvar,hanasusanto2015distributionally,guo2017convergence}. A possible issue with these robust approaches is that they might be too ``pessimistic," especially for tightly constrained setting, which can cause an upwards shift in the optimal value or even feasibility issues; see \cite{RoysetChenEckstrand.24} and examples below. 
Recently, ``optimistic" approaches have emerged to reveal possibilities hidden for their ``pessimistic" counterparts; see, e.g., \cite{nguyen2019optimistic} for likelihood approximation, \cite{jiang2024distributionally} for connections to robust statistics, and \cite{hanasusanto2017ambiguous} for chance-constrained programs. 
Aligning with this ``optimistic" perspective, a series of works \cite{RoysetChenEckstrand.24,antil2024rockafellian,deride2024approximations} mitigate distributional perturbations by solving certain Rockafellian relaxations. 
Rockafellian relaxation involves embedding a minimization problem into a parametric family of problems---an idea pioneered in \cite{Rockafellar.63,Rockafellar.70,Rockafellar.74}; see \cite{Royset.21} and \cite[Chapter 5]{primer} for recent terminology.
Compared with the ``pessimistic'' perspective, these Rockafellian relaxation approaches are more aligned with traditional robust statistics, a diametrically opposing perspective to distributionally robust formulations that is also recognized in \cite{blanchet2024stability,blanchet2024distributionally}.
However, the works \cite{RoysetChenEckstrand.24,antil2024rockafellian,deride2024approximations} do not apply to the chance-constrained setting due to possibly discontinuous integrands, and they also require the concerned distributions to be finite, discrete, or continuous.

In this paper, we report on new constructions within the framework of Rockafellian relaxation and develop new convergence results.
 Our results are applicable to general Borel probability distributions, to lower semicontinuous integrands (with chance-constrained programs as an important application), and to distributional uncertainty quantified by various probability metrics, such as the bounded Lipschitz, Fortet-Mourier, Wasserstein, total variation, and Kullback-Leibler divergence.
A key component of the Rockafellian relaxation framework is the design of approximating Rockafellians, for which we present principled constructions under different assumptions. 
We show that minimizing our approximating Rockafellians recovers the solutions of the actual problem as distributional uncertainty vanishes, without assuming continuity-type conditions like metric regularity.
Another notable feature of our constructions is that, while adopting an ``optimistic'' viewpoint to mitigate distributional inaccuracies, we do not attempt to correct the underlying probability distribution or search for a best distribution over some ambiguity set, which would have led to infinite-dimensional problems. Our formulation is, by design, finite-dimensional, which may offer computational advantages. Moreover, for the chance-constrained setting, we quantify the convergence rate and provide quality guarantees for the approximate solutions. These quantitative stability results are built on assumptions of metric \emph{sub}regularity from variational analysis and upper outer-Minkowski content from geometric measure theory, which are weaker than assumptions in existing works.
As an immediate corollary of our general results, we show that penalized formulations of chance-constrained programs are stable relative to distributional perturbations.

This paper leaves some directions open for future research. We focus on near globally optimal solutions, rather than stationary points, as obtaining similar results for the latter would likely require more stringent assumptions. Nevertheless, our results may already have practical ramifications for mixed-integer stochastic programs in addition to the general theoretical contributions. While we provide a convergence rate for the chance-constrained setting, a general quantitative analysis remains open. This would likely require additional structural assumptions to yield meaningful estimates.

We begin in Section~\ref{sec:general} with general convergence results. Section~\ref{sec:prob-metric} discusses the construction of approximating Rockafellians, considering several types of distributional perturbations. Section~\ref{sec:cc} applies the general results to chance-constrained programs.

\smallskip\smallskip
\noindent {\bf Terminology.} 
For $C\subseteq\reals^n$, we define $\iota_C(x):=0,\bm{1}_C(x):=1$ if $x \in C$; otherwise $\iota_C(x):=\infty,\bm{1}_C(x):=0$. Let $\dist(x, C):=\inf_{y \in C} \|x - y\|_2$ for $x\in\reals^n$. The \emph{boundary} of $C$ is $\bdry C$. Let $\nats:=\{1,2,3,\dots\}$. The \emph{infinite subsets} of $\nats$ is denoted by $\mathcal{N}_\infty^\grill$. For any $\nu \in \nats$, we write $\{1,2,\ldots, \nu\}$ as $[\nu]$. The $m$-dimensional \emph{simplex} is defined as $\Delta_m:=\{p \in \reals^m \mid p \geq 0, \sum_{i=1}^m p_i = 1\}$; also $\Delta_\infty:=\{ p=(p_1, p_2, \ldots) \mid p\geq 0,\sum_{i\in \nats} p_i = 1\}$.  The \emph{Euclidean ball} in $\reals^n$ is $\ball^n(x,\epsilon):=\{y \in \reals^n \mid \|x - y\|_2 \leq \epsilon\}$. For a sequence of sets $C^\nu \subseteq \reals^n$, the \emph{outer limit} is $\nOutLim C^\nu := \{x \in \reals^n \mid \exists N \in \mathcal{N}_\infty^\grill\text{ and }x^\nu \in C^\nu \to x \text{ for } \nu \in N \}$; the \emph{inner limit} is $\nInnLim C^\nu := \{x \in \reals^n \mid \exists x^\nu \in C^\nu \to x\}$; and we say the \emph{set limit} of $C^\nu$ is $C$, denoted by $\nLim C^\nu = C$, when $\nOutLim C^\nu = \nInnLim C^\nu = C$. A set-valued mapping $H:\reals^n \rightrightarrows \reals^d$ is \emph{outer semicontinuous} (osc) if $\bigcup_{x^\nu \to x} \nOutLim H(x^\nu) \subseteq H(x)$ for any $x \in \reals^n$.

A function $f:\reals^n \to \Reals$ is \emph{lower semicontinuous} (lsc) if $\liminf f(x^\nu) \geq f(x)$ for any $x^\nu \to x$; its \emph{effective domain} is $\dom f:=\{x \in \reals^n\mid f(x) < \infty\}$; its \emph{epigraph} is $\epi f := \{(x,t) \in \reals^{n+1}\mid f(x) \leq t\}$; its \emph{lower level sets} are $\lev_{\leq t} f:=\{x \in \reals^n \mid f(x) \leq t\}$; its \emph{minimum value} is $\inf f := \inf \{f(x) \mid x\in \reals^n\}$; and its \emph{sets of minimizers} and \emph{near-minimizers} are $\nargmin f := \{x\in \dom f \mid f(x) = \inf f\}$ and $\epsilon\mbox{-}\nargmin f := \{x\in \dom f \mid f(x) \leq  \inf f+\epsilon\}$, respectively. It is \emph{locally bounded} if for each $x \in \reals^n$, there exist $\epsilon_x>0$ and $M_x < \infty$ such that $|f(y)| \leq M_x$ for all $y \in \ball^n(x,\epsilon_x)$. It is \emph{nondecreasing} if $f(x) \geq f(y)$ for any $x\geq y$ in the element-wise sense. 

The functions $f^\nu:\reals^n \to \Reals$ are said to \emph{epigraphically converge} to a function $f:\reals^n \to \Reals$, denoted by $f^\nu \eto f$, if for every $x \in \reals^n$, one has $\liminf f^\nu(x^\nu) \geq f(x)$ for all $x^\nu\to x$ and $\limsup f^\nu(x^\nu) \leq f(x)$ for some $x^\nu\to x$. The functions $f^\nu : \reals^n \to \Reals$ are \emph{tight} if for any $\epsilon > 0$, there exist a compact set $S_\epsilon \subseteq \reals^n$ and an integer $\nu_\epsilon$ such that $\inf_{S_\epsilon} f^\nu  \leq \inf f^\nu + \epsilon$ for any $\nu \geq \nu_\epsilon$.

Throughout, we assume that $\Xi \subseteq \reals^d$ is nonempty and closed with Borel $\sigma$-algebra $\mathcal{B}(\Xi)$. Let $\mathscr{P}(\Xi)$ be the set of all Borel probability distributions (or simply, distributions) over $\Xi$. 
For $\beta\geq 1$, the set of Borel probability distributions with bounded $\beta$th-order moments is denoted by $\mathscr{P}_\beta(\Xi):=\{\mu \in \mathscr{P}(\Xi) \mid \int_\Xi \|\xi\|^\beta_2\mu(\dd \xi) < \infty\}$. The \emph{support set} of a probability distribution $\mu \in \mathscr{P}(\Xi)$ is denoted by $\supp(\mu)$.  
The \emph{Dirac measure} defined for a point $\xi \in \Xi$ is denoted by $\delta_\xi$.
The \emph{Lebesgue measure} on $\reals^d$ is $\mathcal{L}^d$.
We write $\mu \ll \mathcal{L}^d$ when $\mu$ is \emph{absolutely continuous} with respect to $\mathcal{L}^d$. A function $g:\Xi\times \reals^n\to \Reals$ is \emph{random lsc} if its epigraphical mapping $\xi \mapsto \epi g(\xi, \cdot)$ is measurable in the sense of \cite[Definition 14.1]{VaAn} and closed-valued. For the mappings $G^\nu:\Xi \times \reals^n \to \reals^m$, we say the mappings $G^\nu(\cdot, x)$ are \emph{uniformly bounded} on $\Xi$, \emph{locally uniformly} in $x$, if for any $x \in \reals^n$ there exist $\epsilon_x>0$ and $M_x < \infty$ such that 
\[
\|G^\nu(\xi, y)\|_2 \leq M_x \quad \text{for any $\xi \in \Xi$ and $y \in \ball^n(x,\epsilon_x)$.}
\]

\section{Convergence of Approximating Rockafellians}\label{sec:general}

We recall that $f:\reals^m\times \reals^n\to \Reals$ is a {\em Rockafellian} for $\phi:\reals^n \to \Reals$ if $f(0,x) = \phi(x)$ for all $x\in \reals^n$; see, e.g., \cite[Definition 5.1]{primer}. Throughout, we consider the Rockafellian given by
\begin{equation}\label{eq:f}
	f(u,x) = g_0(x) + h\big(u+\Ex_\mu[ G(\bm{\xi},x) ]\big) + \iota_{\{0\}^m}(u).
\end{equation}
It is immediate that $(u^\star, x^\star) \in \nargmin f$ if and only if $u^\star = 0$ and $x^\star \in \nargmin \phi$. To model changes in the probability distribution $\mu$, we consider distributions $\mu^\nu$ on $\Xi$.
It is easy to construct examples when minimizers and the minimum value of the na\"ive \textit{plug-in functions} $\phi^\nu$, with 
\begin{equation}\label{eq:phi-nu}
\phi^\nu(x) = g_0(x) + h\big(\Ex_{\mu^\nu}[G(\bm{\xi},x) ]\big),
\end{equation}
fail to converge to those of $\phi$, even though $\mu^\nu$ converge to $\mu$ in some sense. 
\begin{example}{\rm (finite distribution I).}\label{ex:finite-I}
	Let $n=m=d=1,  \xi_1 = 0, \xi_2 =1, \Xi=\{\xi_1,\xi_2\}, p_1=p_2=\frac{1}{2}, \mu=p_1\delta_{\xi_1}+p_2\delta_{\xi_2}, g_0(x)=x^2, h=\iota_{(-\infty, 0]}, G(\xi, x)=\frac{1}{2}-\bm{1}_{\{\xi_1\}}(\xi).$ Then, the actual objective function has
	\[
	\phi(x)=x^2+\iota_{(-\infty,0]}\Big(\tfrac{1}{2}-\tfrac{1}{2}\bm{1}_{\{\xi_1\}}(\xi_1)-\tfrac{1}{2}\bm{1}_{\{\xi_1\}}(\xi_2)\Big)=x^2,
	\]
	with $\inf \phi=0$ and $\argmin \phi = \{0\}$.  
	Let $p_1^\nu=p_1-\tfrac{1}{\nu+1}$, $p_2^\nu = p_2+\tfrac{1}{\nu+1}$, and $\mu^\nu=p_1^\nu\delta_{\xi_1}+p_2^\nu\delta_{\xi_2}$, so that $\|p-p^\nu\|_2\to 0$.
	Hence,  the functions $\phi^\nu$ in \eqref{eq:phi-nu} can be written as
	\[
	\phi^\nu(x)=x^2+\iota_{(-\infty,0]}\Big(\tfrac{1}{2}-\big(\tfrac{1}{2}-\tfrac{1}{\nu+1}\big)\bm{1}_{\{\xi_1\}}(\xi_1)-\big(\tfrac{1}{2}+\tfrac{1}{\nu+1}\big)\bm{1}_{\{\xi_1\}}(\xi_2)\Big)=x^2 + \iota_{(-\infty,0]}\big(\tfrac{1}{\nu+1}\big)=\infty.
	\]
	For any $\nu \in \nats$, we have $
	\ninf \phi^\nu = \infty > 0=\ninf \phi$ and $\argmin \phi^\nu = \emptyset$.
\end{example}

To improve stability, we employ an approximation of the Rockafellian $f$ in \eqref{eq:f} given by 
\begin{equation}\label{eq:general-f-nu}
f^\nu(u,x) = g_0(x) + h\big(u + \Ex_{\mu^\nu}[ G^\nu(\bm{\xi},x) ]\big) +  \tfrac{1}{\alpha{\lambda^\nu}}  \|u\|_2^\alpha,
\end{equation}
where $\alpha\geq 1$, $\lambda^\nu > 0$, and $G^\nu:\Xi\times \reals^n\to \reals^m$ with component functions $g^\nu_i:\Xi\times \reals^n\to \reals$ for $i \in [m]$. In fact, $f^\nu$ are Rockafellians of a slight modification of $\phi^\nu$, where $G$ has been replaced by $G^\nu$. We refer to $f^\nu$ as the {\em approximating Rockafellians}. The role of $G^\nu$ and their possible constructions emerge below. 

We establish results about how minimizing $f^\nu$ recovers in the limit as $\nu\to \infty$ solutions of the actual problem \eqref{eq:general-problem}. We establish this under broad conditions in sharp contrast to the situation of for $\phi^\nu$, which requires strong assumptions to be able to recover solutions of \eqref{eq:general-problem} in the limit.

\subsection{Main Results}

We begin with an approximation result for the Rockafellians $f$ and $f^\nu$.
\begin{theorem}{\rm (approximation theorem).}\label{thm:epi-general}
 For  $\mu, \mu^\nu \in \mathscr{P}(\Xi)$, consider a proper, lsc function $g_0:\reals^n \to \Reals$, a proper, lsc, nondecreasing function $h:\reals^m \to \Reals$, and mappings $G,G^\nu:\Xi\times \reals^n \to \reals^m$ with random lsc component functions. Assume that the mappings $G(\cdot, x),G^\nu(\cdot, x)$ are uniformly bounded on $\Xi$, locally uniformly in $x$.
 For fixed $x_0 \in \reals^n$, suppose that the following  conditions hold:
	
	\begin{enumerate}[label=\textnormal{(\roman*)}]
		 \item $G^\nu(\xi,x_0) \leq  G(\xi,x_0)$ for $\mu$-a.e.~$\xi$;\label{item:thm:epi-general-G}
	\item $\liminf \Ex_{\mu^\nu}[ G^\nu(\bm{\xi},x^\nu)] \geq  \Ex_\mu[ G(\bm{\xi},x) ]$ for any $x^\nu \to x$; \label{item:thm:epi-general-liminf}
	 \item $\lambda^\nu \in (0,\infty)\to 0$ and $(\lambda^\nu)^{-1/\alpha}(\Ex_{\mu^\nu}[  G^\nu(\bm{\xi},x_0) ] - \Ex_\mu[ G^\nu(\bm{\xi},x_0) ]) \to 0$. \label{item:thm:epi-general-lambda}
\end{enumerate}
 Then, $f$ and $f^\nu$ in \eqref{eq:f} and \eqref{eq:general-f-nu} are lsc and one has:
\begin{enumerate}[label=\textnormal{(\alph*)}]
	\item $\forall(u^\nu, x^\nu)\to (u,x)$, $\liminf f^\nu(u^\nu,x^\nu) \geq f(u,x)$; \label{item:thm:epi-general-(a)}
	\item $\exists u^\nu \to 0,$  $\limsup f^\nu(u^\nu,x_0)\leq f(0,x_0)$.\label{item:thm:epi-general-(b)}
\end{enumerate}
\end{theorem}

\state Proof. 
Let random lsc $g_{i}, g_{i}^\nu: \Xi \times \reals^n \to \reals$ for $i \in [m]$ be the component functions of $G, G^\nu: \Xi \times \reals^n \to \reals^m$, respectively.
To show $f$ is lsc, it suffices to show that the function $(u,x) \mapsto h(u + \Ex_{\mu}[G(\bm{\xi}, x)])$ is lsc.
Since the measurable functions $g_i(\cdot, x)$ are uniformly bounded on $\Xi$, locally uniformly in $x$, the functions $g_i$ are locally inf-integrable relative to $\mu$; see \cite[p.~541]{primer}.
By \cite[Proposition 8.55]{primer}, the function $x \mapsto \Ex_\mu[g_i(\bm{\xi},x)]$ is lsc for each $i \in [m]$. Let $(u^\nu, x^\nu) \to (u,x)$. By monotonicity of $h$, we conclude that $\liminf h\big(u^\nu + \Ex_{\mu} [G(\bm{\xi},x^\nu)]\big)\geq h\big(u + \Ex_{\mu} [G(\bm{\xi},x)]\big)$. Hence, $f$ is lsc. A similar argument shows that the functions $f^\nu$ are all lsc. We proceed to show Items \ref{item:thm:epi-general-(a)} and \ref{item:thm:epi-general-(b)} as follows.

For Item \ref{item:thm:epi-general-(a)}, consider $(u^\nu, x^\nu) \to (u,x)$. Let $y^\nu=u^\nu + \Ex_{\mu^\nu}[ G^\nu(\bm{\xi},x^\nu)] \in \reals^m$. 
Since the measurable functions $g_i^\nu(\cdot, x)$ are uniformly bounded on $\Xi$, locally uniformly in $x$, we know that the sequence $\{y^\nu\}_\nu$ is bounded. 
For each $i \in [m]$, we have
\begin{equation}\label{prop:epi-general-ineq1}
\infty >\liminf y^\nu_i \geq u_i +  \liminf \Ex_{\mu^\nu} [g_i^\nu(\bm{\xi},x^\nu)] \geq u_i + \Ex_{\mu} [g_i(\bm{\xi},x)] > -\infty,	
\end{equation}
where the third inequality is from \ref{item:thm:epi-general-liminf}.
Let $y,z^\nu\in\reals^m$ be defined by $y_i=\liminf y^\nu_i$ and $z_i^\nu=\inf_{\nu'\geq \nu} y_i^{\nu'}$, so that $z^\nu \to y$ and $z^\nu \leq y^\nu$ for any $\nu \in \nats$. We know that 
\[
\liminf h\big(u^\nu + \Ex_{\mu^\nu} [G^\nu(\bm{\xi},x^\nu)]\big)=\liminf h(y^\nu) \geq \liminf h(z^\nu) \geq h(y) \geq h\big(u + \Ex_{\mu} [G(\bm{\xi},x)]\big) > -\infty,
\]
where the first inequality is by monotonicity of $h$, the second one uses lsc of $h$, the third one is by \eqref{prop:epi-general-ineq1} and monotonicity of $h$, and the last one is from the fact that $h$ is proper.
Similarly, we have  $\liminf g_0(x^\nu) \geq g_0(x) > -\infty$  and $\liminf \frac{1}{\alpha\lambda^\nu}\|u^\nu\|_2^\alpha \geq \iota_{\{0\}^m}(u) > -\infty$. 
Therefore, we conclude that
\[
\begin{aligned}
\liminf f^\nu(u^\nu, x^\nu) &\geq 
\liminf g_0(x^\nu) + \liminf h\big(u^\nu + \Ex_{\mu^\nu} [G^\nu(\bm{\xi},x^\nu)]\big) + \liminf \tfrac{1}{\alpha\lambda^\nu}\|u^\nu\|_2^\alpha \\
&\geq g_0(x) + h\big(u + \Ex_{\mu} [G(\bm{\xi},x)]\big)+\iota_{\{0\}^m}(u) = f(u,x).
\end{aligned}
\]

For Item~\ref{item:thm:epi-general-(b)}, we construct $u^\nu = \Ex_{\mu}[  G^\nu(\bm{\xi},x_0)]-\Ex_{\mu^{\nu}}[ G^\nu(\bm{\xi},x_0) ]$.  From \ref{item:thm:epi-general-lambda}, we have $u^\nu \to 0$ and $\tfrac{1}{\alpha\lambda^\nu}\|u^\nu\|_2^\alpha \to 0$.
Note that 
$
u^\nu+ \Ex_{\mu^{\nu}}[ G^\nu(\bm{\xi},x_0) ] =  \Ex_{\mu}[ G^\nu(\bm{\xi},x_0)].
$
We compute
\[
\begin{aligned}
\limsup f^\nu(u^\nu, x_0) &\leq  g_0(x_0) + \limsup h\big(\Ex_{\mu}[ G^\nu(\bm{\xi},x_0)]\big) + \limsup \tfrac{1}{\alpha\lambda^\nu}\|u^\nu\|_2^\alpha \\
&
\leq g_0(x_0) + h\big(\Ex_{\mu}[ G(\bm{\xi},x_0)]\big)
=  f(0, x_0),  
\end{aligned}
\]
where the second inequality uses \ref{item:thm:epi-general-G} and monotonicity of $h$.
\eop

	Theorem~\ref{thm:epi-general} considers a localized variant of epigraphical convergence, in which the limit supremum part is required to hold only at a fixed point $(0,x_0)$, rather than over all $(u,x)$. Although many of the results and constructions in this paper establish full epigraphical convergence $f^\nu \eto f$, this relaxed one-point approximation may be more suitable in certain applications where full epi-convergence is either nonessential or difficult to guarantee. 
As we see next, the localized version might suffice in practice.

\begin{theorem}{\rm (convergence of approximating Rockafellians).}\label{thm:epi-general-impli}
	Under the conditions of Theorem~\ref{thm:epi-general}, with $x_0 \in \argmin \phi$, the following hold for any $\epsilon \in [0, \infty)$:
	\begin{enumerate}[label=\textnormal{(\alph*)}]
		\item $\ninf f=\ninf \phi$ and $\epsilon\text{-}\nargmin f=\{0\}\times \epsilon\text{-}\nargmin \phi$; \label{item:thm:epi-general-impli-(a)}
		\item $
		\limsup{} ( \ninf  f^\nu) \leq \ninf f
		$;\label{item:thm:epi-general-impli-(b)}
		\item if $\epsilon^\nu \in [0,\infty)\to \epsilon,$ then 
		$
		\nOutLim{} ( \epsilon^\nu\text{-}\nargmin  f^\nu ) \subseteq \epsilon\text{-}\nargmin f
		$;\label{item:thm:epi-general-impli-(c)}

		\item if the sequence $\{(u^\nu,x^\nu) \in \epsilon^\nu\text{-}\nargmin f^\nu\}_\nu$ converges for some $N\in\mathcal{N}_\infty^\grill$ and $\epsilon^\nu \in [0,\infty)\to 0$, then
		$
		\nlim_{\nu \in N}{}  ( \ninf f^\nu ) = \ninf f$; \label{item:thm:epi-general-impli-(d)}
		\item $\ninf f^\nu \to \ninf  f$ if and only if the functions $f^\nu$ are tight; \label{item:thm:epi-general-impli-(e)}
		\item if $\ninf f^\nu \to \ninf f$, then there exist $\epsilon^\nu \in [0,\infty)\to 0$ such that
		$(0,x_0)\in \nInnLim{} (\epsilon^\nu\text{-}\nargmin f^\nu)
		$. \label{item:thm:epi-general-impli-(f)}
	\end{enumerate}
	Moreover, if conditions \ref{item:thm:epi-general-G} and \ref{item:thm:epi-general-lambda} of Theorem~\ref{thm:epi-general} hold for every $x_0 \in \dom g_0$ with the same parameters $\lambda^\nu$, then $f^\nu \eto f$, and consequently, Item~{\rm \ref{item:thm:epi-general-impli-(f)}} can be strengthened as follows:
	\begin{enumerate}[label=\textnormal{(g)}]
		\item if $ \ninf f^\nu \to \ninf f$, then there exist $\epsilon^\nu \in [0,\infty)\to 0$ such that
		$
		\nLim{} ( \epsilon^\nu\text{-}\nargmin f^\nu) = \nargmin f.
		$\label{item:thm:epi-general-impli-(g)}
	\end{enumerate}
\end{theorem}
 
\state Proof. 
Item~\ref{item:thm:epi-general-impli-(a)} follows from the definition of~$f$. By assumption, $\inf f = \inf \phi = f(0, x_0) > -\infty$. The arguments for Items~\ref{item:thm:epi-general-impli-(b)}--\ref{item:thm:epi-general-impli-(f)} are slight modifications of the proof of~\cite[Theorem~5.5]{primer}, combined with the conclusion of Theorem~\ref{thm:epi-general}. Specifically, we apply Theorem~\ref{thm:epi-general}\ref{item:thm:epi-general-(b)} (cf.~\cite[Theorem~4.15(b)]{primer}) at the point $(0, x_0) \in \argmin f$, whose existence is guaranteed by the assumption. For Item~\ref{item:thm:epi-general-impli-(g)}, we observe that $\dom f \subseteq \{0\} \times \dom g_0$, so that $f^\nu \eto f$ by Theorem~\ref{thm:epi-general} and  \cite[Proposition 7.7]{VaAn} applies.
\eop

\subsection{Discussion}\label{sec:discussion}
The approximating Rockafellians $f^\nu$ in \eqref{eq:general-f-nu} involve three distinct forms of approximation: the distributions $\mu^\nu$ to $\mu$, the mappings $G^\nu$ to $G$, and the functions $\tfrac{1}{\alpha\lambda^\nu}\|\cdot\|_2^\alpha$ to $\iota_{\{0\}^m}$. 
As highlighted in Section~\ref{sec:intro}, the central goal of this paper is to stabilize the problem under possibly adversarial perturbations to the actual  distribution $\mu$. Guided by the variational convergence in Theorem~\ref{thm:epi-general-impli}, our approach is to \textit{design} the parameters $\lambda^\nu$ and mappings $G^\nu$ so as to actively counteract inaccuracies introduced by $\mu^\nu$, which we will detail in Sections~\ref{sec:prob-metric} and \ref{sec:cc}. But before that, 
in this subsection, we consider several illustrative examples in which the distributions $\mu$ and $\mu^\nu$ are supported on finite or discrete sets. These examples allow us to concretely demonstrate the mitigating effects afforded by approximating Rockafellians $f^\nu$ and also facilitate comparisons with existing approaches in the literature.

\subsubsection{Finite Distribution}

To begin with, assume that the distribution $\mu$ is supported on a finite set $\Xi=\{\xi_k \in \reals^d \mid k \in [s]\}$ and can be represented as $\mu = \sum_{k=1}^s p_k\delta_{\xi_k}$ with $p \in \Delta_s$. This is a setting similar to the one studied in \cite[Section 4.3]{RoysetChenEckstrand.24}. 
When $\mu$ is approximated by $\mu^\nu = \sum_{k=1}^s p_k^\nu\delta_{\xi_k}$ with $p^\nu \in \Delta_s$, the plug-in function can be written as
$\phi^\nu(x)=g_0(x) + h(\nsum_{k=1}^s p_k^\nu G(\xi_k, x)).$
As shown in Example~\ref{ex:finite-I},
even for finite distributions $\mu$ and $\mu^\nu$ with $p^\nu \to p$, 
the solution sets $\argmin \phi^\nu$ can be always empty. 
Complementing Example~\ref{ex:finite-I}, the following example highlights a more subtle failure mode: even when $\inf \phi^\nu < \infty$ and $\argmin \phi^\nu \neq \emptyset$ for all $\nu \in \nats$, neither  $\inf \phi^\nu$ nor  $\argmin \phi^\nu$ necessarily converge to those of the actual problem. This illustrates that minimizing $\phi^\nu$ may yield erroneous estimations that are difficult to diagnose.

\begin{example}{\rm (finite distribution II).}\label{ex:finite-II}
Consider the same setting as in Example~\ref{ex:finite-I}, except with $G(\xi, x)=\frac{1}{2}-\bm{1}_{\{\xi_1+x\}}(\xi).$ We can write the actual $\phi$ in \eqref{eq:general-problem} as
	\[
	\phi(x)=x^2+\iota_{(-\infty,0]}\Big(\tfrac{1}{2}-\tfrac{1}{2}\bm{1}_{\{\xi_1+x\}}(\xi_1)-\tfrac{1}{2}\bm{1}_{\{\xi_1+x\}}(\xi_2)\Big)=x^2 + \iota_{\{0,\xi_2-\xi_1\}}(x).
	\]
	Then, $\inf \phi=0$ and $\argmin \phi = \{0\}$.  
	For $\nu \in \nats$, let $p_1^\nu=p_1-\tfrac{1}{\nu+1}$ and $p_2^\nu = p_2+\tfrac{1}{\nu+1}$, so that $\|p-p^\nu\|_2\to 0$.
	Hence, the plug-in functions $\phi^\nu$ in \eqref{eq:phi-nu} can be written as
	\[
	\phi^\nu(x)=x^2+\iota_{(-\infty,0]}\Big(\tfrac{1}{2}-\big(\tfrac{1}{2}-\tfrac{1}{\nu+1}\big)\bm{1}_{\{\xi_1+x\}}(\xi_1)-\big(\tfrac{1}{2}+\tfrac{1}{\nu+1}\big)\bm{1}_{\{\xi_1+x\}}(\xi_2)\Big)=x^2 + \iota_{\{\xi_2-\xi_1\}}(x).
	\]
	For any $\nu \in \nats$, we have $\ninf \phi^\nu = 1 > 0=\ninf \phi$ and $\argmin \phi^\nu = \{1\} \not\subseteq\{0\}=\argmin \phi$. 
\end{example}

One possible remedy for these disproportionate changes in the solutions is to minimize, instead of $\phi^\nu$ in \eqref{eq:phi-nu}, the approximating Rockafellians $f^\nu$ defined in \eqref{eq:general-f-nu}, 
which resembles the construction in \cite[(4.16)]{RoysetChenEckstrand.24}, except that the mapping $G$ here is permitted to be discontinuous. As a corollary of Theorem~\ref{thm:epi-general-impli} extending \cite[Proposition 4.8]{RoysetChenEckstrand.24} to discontinuous integrands, we obtain the following result.
\begin{corollary}{\rm (finite distribution).}\label{coro:finite}
	Suppose that $\Xi = \{\xi_k \mid k \in [s]\}$. Let $g_0: \reals^n \to \Reals$ be a proper, lsc function, $h: \reals^m \to \Reals$ be a proper, lsc, and nondecreasing function, and $G$ be a mapping from $\Xi\times \reals^n$ to $\reals^m$. Assume that the component functions $g_i(\xi_k, \cdot)$ of $G(\xi_k,\cdot)$ are locally bounded and lsc for each $k \in [s]$. For $f$ and $f^\nu$ in \eqref{eq:f} and \eqref{eq:general-f-nu} with $G^\nu=G$, if $\argmin \phi \neq \emptyset$ for $\phi$  in \eqref{eq:general-problem}, $\lambda^\nu \in (0,\infty) \to 0$ and $\tfrac{1}{\lambda^\nu}\|p - p^\nu\|_2^\alpha \to 0$ with $p,p^\nu \in \Delta_s$, then Theorem~\ref{thm:epi-general-impli}\ref{item:thm:epi-general-impli-(a)}--\ref{item:thm:epi-general-impli-(g)} hold and $f^\nu \eto f$.	
\end{corollary}

\state Proof.
Let $\mu = \sum_{k=1}^s p_k\delta_{\xi_k}$, $\mu^{\nu} = \sum_{k=1}^s p_k^\nu\delta_{\xi_k}$, and $G^\nu = G$.
Note that the mapping $G:\Xi \times \reals^n \to \reals^m$ is uniformly bounded on $\Xi$, locally uniformly in $x$, if and only if $g_i(\xi_k, \cdot)$ is locally bounded for any $k \in [s]$ and $i\in[m]$.
We only need to verify that the conditions in Theorem~\ref{thm:epi-general} hold for any $x \in \dom g_0$. 
Since $G^\nu = G$, Theorem~\ref{thm:epi-general}\ref{item:thm:epi-general-G} holds trivially. For Theorem~\ref{thm:epi-general}\ref{item:thm:epi-general-liminf}, we have
\[
\liminf \Ex_{\mu^\nu}[g_i(\bm{\xi},x^\nu)] = \nsum_{k=1}^s p_k \liminf  g_i(\xi_k,x^\nu)\geq 
\nsum_{k=1}^s p_k g_i(\xi_k,x)
=\Ex_{\mu}[g_i(\bm{\xi},x)],
\] 
where the inequality follows from the fact that $g_i$ is random lsc. Note that
\[
|\Ex_{\mu^\nu}[g_i(\bm{\xi},x)] - \Ex_{\mu}[g_i(\bm{\xi},x)]|=\Big| \nsum_{k=1}^s (p_k^\nu-p_k) g_i(\xi_k,x) \Big| \leq \sqrt{s}\max_{k\in [s]}|g_i(\xi_k,x)|\|p - p^\nu\|_2.
\]
If $\lambda^\nu \in (0,\infty) \to 0$ and $\tfrac{1}{\lambda^\nu}\|p - p^\nu\|_2^\alpha \to 0$, then $(\lambda^\nu)^{-1/\alpha}(\Ex_{\mu^\nu}[G(\bm{\xi},x)] - \Ex_{\mu}[G(\bm{\xi},x)])\to 0$ and hence Theorem~\ref{thm:epi-general}\ref{item:thm:epi-general-lambda} holds. The conclusion then follows by invoking Theorem~\ref{thm:epi-general-impli}.
\eop

\subsubsection{Discrete Distribution}
We now consider a more general setting in which $\mu$ is supported on a countable set $\Xi=\{\xi_k \in \reals^d\mid k \in \nats\}$ and admits the representation $\mu = \sum_{k=1}^\infty p_k \delta_{\xi_k}$ with $p=(p_1,p_2,\ldots) \in \Delta_\infty$. To approximate $\mu \in \mathscr{P}(\Xi)$, we introduce $\mu^\nu = \sum_{k=1}^\infty p_k^\nu \delta_{\xi_k}$ with $p^\nu = (p_1^\nu,p_2^\nu,\ldots) \in \Delta_\infty$ for $\nu \in \nats$. 
This countable setting allows us to examine more general modes of convergence of $\mu^\nu$ to $\mu$. 

For comparison with the finite distribution setting, especially Corollary~\ref{coro:finite}, we temporarily restrict ourselves to using the approximating functions $\tfrac{1}{\alpha\lambda^\nu}\|\cdot\|_2^\alpha$ only, without exploiting the flexibility of designing $G^\nu$. The following example demonstrates that when $\mu^\nu$ converge weakly to $\mu$ (rather than $p^\nu \to p$ in a suitable norm), unlike the conclusion of Corollary~\ref{coro:finite}, for any $\lambda^\nu \in (0,\infty) \to 0$, we have neither $\limsup{} (\inf f^\nu) \leq \inf f$ nor $\nOutLim{} (\argmin f^\nu) \subseteq \argmin f$. 
We also illustrate the mitigating effect of appropriately designed $G^\nu$ on these disproportionate sensitivities.

\begin{example}{\rm (discrete distribution I).}\label{ex:discreteI}
	Consider the following setting: $n=m=d=1,  \xi_1=1, \xi_k = 1+\tfrac{1}{k}$ for $k \geq 2$, $\Xi=\{\xi_k\mid k \in \nats\}, p_1=1, p_k=0$ for $k \geq 2$, $g_0(x)=x^2, h(u)=\iota_{(-\infty, 0]}(u), G(\xi, x)=\frac{1}{2}-\bm{1}_{[0,\xi_1]\times [1,2]}(\xi,x).$ Let $\mu=\sum_{k=1}^\infty p_k\delta_{\xi_k}$. Since $G(\xi,x)=\tfrac{1}{2}$ whenever $x \notin [1,2]$, we can write 
	\[
	\phi(x)=x^2+\iota_{(-\infty,0]}\big(\tfrac{1}{2}-\mu(\{\xi\in \Xi\mid\xi\in[0,\xi_1],x\in[1,2]\})\big)=x^2 + \iota_{[1,2]}(x).
	\]
	Then, $\inf \phi=1$ and $\argmin \phi = \{1\}$.  
	For $\nu\in\nats$, let $p_k^\nu = \bm{1}_{\{k\}}(\nu)$ for any $k \in \nats$ and $\mu^\nu=\sum_{k=1}^\infty p_k^\nu \delta_{\xi_k}$, so that $\mu^\nu$ converge to $\mu$ weakly
	but not in total variation.
	Hence, for any $\nu  \geq 2$, the approximating Rockafellians $f^\nu$ in \eqref{eq:general-f-nu} with $G^\nu=G$ can be written as
	\[
	\begin{aligned}
	f^\nu(u,x)&=x^2+\iota_{(-\infty,0]}\big(u+\tfrac{1}{2}-\mu^\nu(\{\xi\in \Xi\mid\xi\in[0,\xi_1],x\in[1,2]\})\big)+\tfrac{1}{\alpha\lambda^\nu}|u|^\alpha\\
	&=x^2+\iota_{(-\infty,-1/2]}(u)+\tfrac{1}{\alpha\lambda^\nu}|u|^\alpha.
	\end{aligned}
	\]
	For any $\nu  \geq 2$, a simple computation reveals that 
	$\inf f^\nu = (\alpha 2^\alpha \lambda^\nu)^{-1}$ and $\argmin f^\nu=\{(-\tfrac{1}{2},0)\}$.
	Hence, for any $\alpha \geq 1$ and $\lambda^\nu \in (0,\infty)\to 0$, we have
		\[
	\nlimsup{} (\ninf f^\nu)=\infty > 1= \ninf \phi,\quad
	\nOutLim{} (\argmin f^\nu)=\{(-\tfrac{1}{2},0)\}\not\subseteq
	\{(0,1)\}=\argmin f.
	\]
	These instabilities can be addressed by appropriately designing approximating mappings $G^\nu$.
\end{example}
\state Detail.
A possible construction of $G^\nu:\Xi \times \reals^n \to \reals^m$ is $G^\nu(\xi, x)=\frac{1}{2}-\bm{1}_{[0,\xi_1+\nu^{-1}]\times [1,2]}(\xi,x)$. Note that $G^\nu(\xi, x) \leq G(\xi, x)$ for any $\nu \in \nats$. For any $x^\nu \to x$, we have
$
\liminf \Ex_{\mu^\nu}[G^\nu(\bm{\xi}, x^\nu)] =\liminf{} \tfrac{1}{2}-\bm{1}_{ [1,2]}(x^\nu)\geq\Ex_\mu[G(\bm{\xi}, x)].
$
Finally, we know that 
$
\Ex_{\mu^\nu}[G^\nu(\bm{\xi}, x)] - \Ex_{\mu}[G^\nu(\bm{\xi}, x)] = 0
$ always. Therefore, the convergence can be obtained by invoking Theorem~\ref{thm:epi-general-impli} with any $\lambda^\nu \in (0,\infty)\to 0$.
\eop

Similar to Example~\ref{ex:finite-II} in the finite distribution setting, the following example illustrates that, in the absence of properly designed mappings $G^\nu$, minimizing $f^\nu$ may also yield erroneous estimates of minimum values and solutions that are difficult to detect.

\begin{example}{\rm (discrete distribution II).}\label{ex:discreteII}
	Consider the same setting as in Example~\ref{ex:discreteI}, except with $g_0(x)=-\mathbf{1}_{(-\infty,\xi_1]}(x) - \tfrac{1}{2}\mathbf{1}_{[\xi_2,\infty)}(x)$ and $G(\xi, x)=\frac{1}{2}-\bm{1}_{(-\infty,x]}(\xi)$. We can write $\phi$ as
	\[
	\begin{aligned}
	\phi(x) &= -\mathbf{1}_{(-\infty,\xi_1]}(x) - \tfrac{1}{2}\mathbf{1}_{[\xi_2,\infty)}(x) + \iota_{(-\infty, 0]}\big(\tfrac{1}{2}-\mu(\{\xi \in \Xi \mid x\geq \xi\})\big) \\
	&=-\mathbf{1}_{(-\infty,\xi_1]}(x) - \tfrac{1}{2}\mathbf{1}_{[\xi_2,\infty)}(x) + \iota_{[\xi_1,\infty)}(x).
	\end{aligned}
	\]
	Then, $\inf \phi=-1$ and $\argmin \phi = \{\xi_1\}$.  
	For $\nu\in\nats$, let $p_k^\nu = \bm{1}_{\{k\}}(\nu)$ for each $k \in \nats$ and $\mu^\nu=\sum_{k=1}^\infty p_k^\nu \delta_{\xi_k}$, so that $\mu^\nu$ converge to $\mu$ weakly. 
	Hence,  $f^\nu$ in \eqref{eq:general-f-nu} with $G^\nu=G$ can be written as
	\[
	f^\nu(u,x)=-\mathbf{1}_{(-\infty,\xi_1]}(x)- \tfrac{1}{2}\mathbf{1}_{[\xi_2,\infty)}(x)+\iota_{(-\infty,0]}\big(u+\tfrac{1}{2}-\mu^\nu(\{\xi\in\Xi\mid x\geq \xi\})\big)+\tfrac{1}{\alpha\lambda^\nu}|u|^\alpha.
	\]
	Whenever $\lambda^\nu \in (0,(\alpha2^\alpha)^{-1})$, we observe that $\inf f^\nu = -1/2$ with $\argmin f^\nu = \{0\}\times [\xi_2, \infty)$.
	Hence, for any $\alpha \geq 1$ and $\lambda^\nu \in (0,\infty) \to 0$, we have
	\[
	\nlimsup{} (\ninf f^\nu)=-1/2>-1= \ninf \phi,\quad
	\nOutLim{} (\argmin f^\nu)=\{0\}\times[\xi_2, \infty)\not\subseteq
	\{(0,\xi_1)\}=\argmin f.
	\]
	These instabilities can be addressed by appropriately designing approximating mappings $G^\nu$.
\end{example}
\state Detail. A possible construction of $G^\nu:\Xi\times \reals^n\to \reals^m$ is $G^\nu(\xi, x)=\frac{1}{2}-\bm{1}_{(-\infty, x+\nu^{-1}]}(\xi)$. Note that $G^\nu(\xi, x) \leq G(\xi, x)$ for any $\nu \in \nats$. 
For $\nu  \geq 2$, we observe that
$
\Ex_{\mu^\nu}[G^\nu(\bm{\xi}, x^\nu)]-\tfrac{1}{2}=-\bm{1}_{(-\infty, x^\nu+\nu^{-1}]}(\xi_1+\nu^{-1})=-\bm{1}_{[\xi_1,\infty)}(x^\nu).
$
Hence, by lsc of the function $x\mapsto -\bm{1}_{[\xi_1,\infty)}(x)$, we have
$
\liminf \Ex_{\mu^\nu}[G^\nu(\bm{\xi}, x^\nu)]\geq \tfrac{1}{2}-\bm{1}_{[\xi_1,\infty)}(x)= \Ex_\mu[G(\bm{\xi}, x)].
$
Finally, for sufficiently large $\nu \in \nats$, we know that
\[
\Ex_{\mu^\nu}[G^\nu(\bm{\xi}, x)] - \Ex_{\mu}[G^\nu(\bm{\xi}, x)] = \bm{1}_{(-\infty, x]}(\xi_1-\nu^{-1}) - \bm{1}_{(-\infty, x]}(\xi_1)=0.
\]
Therefore, the convergence can be obtained by invoking Theorem~\ref{thm:epi-general-impli} with any $\lambda^\nu \in (0,\infty)\to 0$.
\eop

\begin{remark}{\rm (alternative Rockafellian).}\label{rmk:alternative}
	The discrete distribution setting is also studied in \cite[Section 4.4]{RoysetChenEckstrand.24} for a slightly different problem without the outer function $h$. Motivated by \cite[Section 4.4]{RoysetChenEckstrand.24}, we consider an alternative Rockafellian and its approximations as follows: 
	\begin{align}
	f(u,x)&=g_0(x)+h\Big(\nsum_{k=1}^\infty (p_k^\nu+u_k)G(\xi_k,x)\Big) + \iota_{\{0\}^m}(u),\\
	f^\nu(u,x)&=g_0(x)+h\Big(\nsum_{k=1}^\infty (p_k^\nu+u_k)G(\xi_k,x)\Big) + \tfrac{1}{\lambda^\nu}\|u\|_{\ell_1}+\iota_{\Delta_\infty}(p^\nu + u).\label{eq:alternative-f-nu}
\end{align}
When $\mu^\nu$ converge to $\mu$ only weakly and using similar constructions as in Examples~\ref{ex:discreteI} and \ref{ex:discreteII}, one can readily verify  that for any $\lambda^\nu \in (0,\infty)\to 0$ and $f^\nu$ in \eqref{eq:alternative-f-nu}, neither $\inf f^\nu$ nor  $\argmin f^\nu$ necessarily converge to those of the Rockafellian $f$. As will be shown later, this deviation from \cite[Proposition 4.9]{RoysetChenEckstrand.24} arises because $\mu^\nu$ do not converge to $\mu$ in $\ell_1$-norm, since $\|p^\nu - p\|_{\ell_1} = 2$ always for large $\nu \in \nats$.
\end{remark}

In sum, to accommodate a broad class of distributions $\mu^\nu$ that converge to $\mu$ with varying strength, it is beneficial to construct approximating functions $\tfrac{1}{\alpha\lambda^\nu}\|\cdot\|_2^\alpha$ and mappings $G^\nu$ that mesh appropriately with $\mu^\nu$ in the sense of Theorem~\ref{thm:epi-general-impli}. In the next section, we illustrate several principled constructions for designing $G^\nu$ and selecting $\lambda^\nu$ under different modes of convergence of probability distributions.
 
\section{Convergence of Probability Distributions}\label{sec:prob-metric}

To apply Theorem~\ref{thm:epi-general-impli}, the mappings $G^\nu$ and parameters $\lambda^\nu$ must satisfy the conditions described in Theorem~\ref{thm:epi-general}\ref{item:thm:epi-general-G}--\ref{item:thm:epi-general-lambda}. 
In this section, we present several principled constructions
assuming various modes of convergence of $\mu^\nu$ to $\mu$. For notational simplicity, we confine the development to a single component function of the mapping $G$. Let $g:\Xi \times \reals^n \to \reals$ be a random lsc function. 
The distributions $\mu^\nu$ are assumed to converge to $\mu$ in some sense specified later. For a fixed $x_0 \in \reals^n$, the main goal of this section is then to construct $g^\nu : \Xi \times \reals^n \to \reals$ satisfying: 
\begin{align}
&	g^\nu(\xi,x_0) \leq g(\xi, x_0) \mbox{ for } \mu\mbox{-a.s}.~\xi;\label{eq:general-conditions-i}\\ 
&	\liminf \Ex_{\mu^\nu}[g^\nu(\bm{\xi},x^\nu)] \geq \Ex_\mu[g(\bm{\xi},x)] \mbox{ for any } x^\nu \to x;\label{eq:general-conditions-ii}\\
&	\Ex_{\mu^\nu}[g^\nu(\bm{\xi},x_0)] - \Ex_\mu[g^\nu(\bm{\xi},x_0)] \to 0.\label{eq:general-conditions-iii}
\end{align}
These conditions align with  Theorem~\ref{thm:epi-general}\ref{item:thm:epi-general-G}--\ref{item:thm:epi-general-lambda}. We begin with existence results, assuming $\mu^\nu$ converging in the weak or setwise sense.  
Then, we provide explicit guidelines by analyzing different probability metrics. We conclude this section with a discussion of the empirical approximation.

\subsection{Weak Convergence}\label{sec:weak}

Recall that $\mu^\nu \in \mathscr{P}(\Xi)$ converge to $\mu \in \mathscr{P}(\Xi)$ weakly if and only if $\Ex_{\mu^\nu}[\ell(\bm{\xi})] \to \Ex_{\mu}[\ell(\bm{\xi})]$ for every bounded and continuous function $\ell:\Xi \to \reals$. 
The inequality in \eqref{eq:general-conditions-i} is usually easy to verify.
The following proposition presents sufficient conditions validating \eqref{eq:general-conditions-ii} and \eqref{eq:general-conditions-iii}.

\begin{proposition}\label{prop:g-perturb}
 Suppose that $\mu^\nu \in \mathscr{P}(\Xi)$ converge to $\mu \in \mathscr{P}(\Xi)$ weakly.  For functions $g(\cdot, x), g^\nu(\cdot, x)$ that are uniformly bounded on $\Xi$, locally uniformly in $x$, we assume the following conditions hold:
\begin{enumerate}[label=\textnormal{(\roman*)}]
	\item the function $\xi\mapsto g^\nu(\xi, x)$ is continuous on $\Xi$ for any $x$ and $\nu\in \nats$;\label{item:lem:g-perturb-cont} 
	\item for any $\xi^\nu \to \xi$ in $\Xi$ and $x^\nu\to x$, $\nliminf g^\nu(\xi^\nu,x^\nu) \geq g(\xi, x)$.\label{item:lem:g-perturb-liminf}
	\end{enumerate}	
Then, one has: 
\begin{enumerate}[label=\textnormal{(\alph*)}]
	\item for any $N \in  \mathcal{N}_\infty^\grill$ and $\lim_{\nu \in N} x^\nu = x$, $\liminf_{\nu \in N} \Ex_{\mu^\nu}[g^\nu(\bm{\xi},x^\nu) ] \geq \Ex_\mu[ g(\bm{\xi},x) ]$; \label{item:prop:g-perturb-(a)}
\item for any fixed $x_0$, there exist $i_\nu^\circ\in \nats\uparrow \infty$ such that for any $i_\nu \in \nats\uparrow \infty$ satisfying $i_\nu \leq i_\nu^\circ$ and sufficiently large $\nu$, we have
$
|\Ex_{\mu^\nu}[ g^{i_\nu}(\bm{\xi}, x_0) ] - \Ex_\mu[ g^{i_\nu}(\bm{\xi}, x_0) ]| \leq \tfrac{1}{i_\nu}\to 0.
$ \label{item:prop:g-perturb-(b)}
\end{enumerate}	
\end{proposition}
\state Proof. 
For Item~\ref{item:prop:g-perturb-(a)} and fixed $N \in  \mathcal{N}_\infty^\grill$, since $\{x^\nu\}_{\nu \in N}$ can be naturally extended to $\{\bar{x}^\nu\}_{\nu\in \nats}$ such that $\bar{x}^\nu \to x$, noting $\liminf_{\nu \in N} \Ex_{\mu^\nu}[g^\nu(\bm{\xi},x^\nu) ] \geq \liminf_{\nu \in \nats} \Ex_{\mu^\nu}[g^\nu(\bm{\xi},\bar{x}^\nu) ] $, we can simply assume $N=\nats$.
From \ref{item:lem:g-perturb-cont}, we know that the functions $g^\nu(\cdot,x)$ are Borel measurable for any $x$. Hence, using the boundedness and an extended Fautou's lemma for weakly convergent measures \cite[Theorem 3.4]{feinberg2022epi}, we have
\[
\liminf \Ex_{\mu^\nu}[  g^\nu(\bm{\xi},x^\nu)] \geq \Ex_\mu\Big[ \liminf_{\nu \in \nats,\xi^\nu\to \bm{\xi}} g^\nu(\xi^\nu,x^\nu)  \Big]\geq\Ex_\mu[ g(\bm{\xi},x) ],
\]
where the last inequality is from \ref{item:lem:g-perturb-liminf}.

For Item~\ref{item:prop:g-perturb-(b)}, as $\mu^\nu$ converge to $\mu$ weakly, and the functions $g^\nu(\cdot, x_0)$ are bounded and continuous on $\Xi$,  we have $\Ex_{\mu^\nu}[g^{a}(\bm{\xi}, x_0)] \to \Ex_{\mu}[g^a(\bm{\xi}, x_0)]$ as $\nu \to \infty$ for any fixed $a \in \nats$.  Define
\[
\eta_{a,b}=\inf\bigg\{ 
\nu \in \nats \biggm| \sup_{\nu' \geq \nu}\Big| \int_\Xi g^a(\xi, x_0) (\mu^{\nu'} - \mu)(\dd \xi) \Big| \leq \frac{1}{b}
\bigg\}.
\]
For every $a,b \in \nats$, we know that $\eta_{a,b} < \infty$. 
Let us also define
$
\eta^k =\max \big\{ \eta_{k',k} \bigm| k'\in [ k]\big\},
$
which is finite for any $k \in \nats$.
We proceed by discussing two cases.

If the sequence $\{\eta^k\}_k$ is bounded above, say, by $\eta \in \nats$, then
we have $|\int_\Xi g^{k'}(\xi, x_0) (\mu^\nu-\mu)(\dd \xi)| \leq \tfrac{1}{k}$ for any $\nu \geq \eta$, any $k \in \nats$, and any $k' \leq k$.  
Therefore, for $\nu \geq \eta$, it holds that
$
| \int_\Xi g^k(\xi, x_0) (\mu^\nu - \mu)(\dd \xi) | 
=0
$
uniformly for all $k \in \nats$. Hence, we can choose arbitrary $i_\nu^\circ\uparrow \infty$.

Suppose that $\{\eta^k\}_k$ is unbounded; i.e., for any $\nu$, there exists $k \in \nats$ such that $\eta^k \geq \nu$.  Let 
\[
i_\nu^\circ=
\left\{ \begin{array}{rcl}
         1 & \mbox{for}
         & \nu < \eta^1, \\
         \max\{k \mid \nu \geq \eta^{k}\} & \mbox{for} & \nu \geq \eta^1.
                \end{array}\right.
\]
Since $\{\eta^k\}_{k}$ is unbounded, we can see that $i_\nu^\circ < \infty$ for all $\nu \in \nats$ and $i_\nu^\circ \uparrow \infty$ as $\nu \to \infty$. Moreover, 
for any $\nu \geq \eta^1$, $i_\nu \leq i_\nu^\circ$, and $i_\nu \uparrow \infty$, we have
$
| \int_\Xi g^{i_\nu}(\xi, x_0) (\mu^\nu - \mu)(\dd \xi) |\leq \frac{1}{i_\nu^\circ} \leq \frac{1}{i_\nu} \to 0
$
as $\nu \to \infty$.
\eop

The continuous functions $g^\nu$ might be expressed using parameters that control the fidelity of the approximation.  Proposition~\ref{prop:g-perturb} shows that, under mild assumptions, there always exists a sequence of parameters such that \eqref{eq:general-conditions-ii} and \eqref{eq:general-conditions-iii} hold. If the functions $g^\nu$ also satisfy  \eqref{eq:general-conditions-i}, and $\lambda^\nu$ are chosen according to \eqref{eq:general-conditions-iii}, then conditions Theorem~\ref{thm:epi-general}\ref{item:thm:epi-general-G}--\ref{item:thm:epi-general-lambda} can be verified similarly, and convergence of solutions can be guaranteed by invoking Theorem~\ref{thm:epi-general-impli}.

If the function $g$ is merely random lsc, then it appears difficult, in general, to design approximating functions $g^\nu$. Thus, the next subsection assumes that $(\xi,x)\mapsto g(\xi,x)$ is lsc jointly in both $\xi$ and $x$. This leads to a principled construction based on epigraphical regularization.

\subsubsection{Epigraphical  Regularization}\label{sec:envelope}

For an lsc function $g:\Xi \times \reals^n \to \reals$ with $g(\cdot, x)$ uniformly bounded on $\Xi$, locally uniformly in $x$, consider the following functions $g^\nu$ with parameter $\theta^\nu \in (0,\infty) \to 0$ and exponent $\beta \geq 1$: \begin{equation}\label{eq:pth-envelope}
g^\nu(\xi, x)=\ninf_{\zeta \in \Xi}  g(\zeta,x)+\tfrac{1}{\beta\theta^\nu} \|\xi - \zeta\|_2^\beta.
\end{equation}	
For each $x \in \reals^n$, the function $g(\cdot, x)$ is bounded on $\Xi$, and hence $g^\nu(\cdot, x) > -\infty$.

We refer to $g^\nu$ in \eqref{eq:pth-envelope} as the \emph{epigraphical regularization} of $g$ with parameters $\beta$ and $\theta^\nu$, which can be viewed as the epigraphical sum (inf-convolution) of the function $g$ and the kernel $(\xi,x)\mapsto \tfrac{1}{\beta\theta^\nu}\|\xi\|_2^\beta$; see \cite[Section 3]{AttouchWets.91}. 
Computationally, it may be helpful to have the flexibility of choosing different $\beta \geq 1$ in  \eqref{eq:pth-envelope}.
Two notable special cases, corresponding to $\beta = 1$ and $\beta = 2$, are described below.

\begin{example}{\rm (Pasch-Hausdorff partial envelope).}\label{example:PH}
	Let $\beta=1$ and $x \in \reals^n$. We can write \eqref{eq:pth-envelope} as  
\[
g^\nu(\xi, x)=\ninf_{\zeta \in \Xi} g(\zeta,x)+\tfrac{1}{\theta^\nu} \|\xi - \zeta\|_2,
\]
which is known as the Pasch-Hausdorff envelope \cite[Example 9.11]{VaAn} of $\xi \mapsto g(\xi, x) + \iota_\Xi(\xi)$. This is also related to the Baire-Wijsman approximation and Lipschitz regularization in the literature. 
\end{example}

\begin{example}{\rm (Moreau partial envelope).}\label{example:moreau}
Let $\beta=2$ and $x \in \reals^n$. We can write \eqref{eq:pth-envelope} as  
\[
g^\nu(\xi, x)=\ninf_{\zeta \in \Xi} g(\zeta,x)+\tfrac{1}{2\theta^\nu} \|\xi - \zeta\|_2^2,
\]
which is known as the Moreau envelope \cite[Definition 1.22]{VaAn} of  $\xi \mapsto g(\xi, x) + \iota_\Xi(\xi)$; see also \cite[Section 3.2]{li2024decomposition} for computational usage of the Moreau partial envelope.
\end{example}

We now collect some properties of the epigraphical regularization \eqref{eq:pth-envelope}.

\begin{proposition}\label{prop:p-env}
Let $g:\Xi\times \reals^n\to \reals$ be lsc. Assume that for each $x \in \reals^n$, there exist $\epsilon_x > 0$ and $M_x < \infty$ such that $|g(\xi, y)| \leq M_x$ for any $\xi\in \Xi$ and $y\in \ball^n(x,\epsilon_x)$. Then, for the functions $g^\nu : \Xi \times \reals^n \to \reals$ in \eqref{eq:pth-envelope} with $\theta^\nu\in(0,\infty)\to 0$, the following hold:
\begin{enumerate}[label=\textnormal{(\alph*)}]
		\item the functions $g^\nu $ are lsc with $|g^\nu(\xi, y)| \leq M_x$ for any $x \in \reals^n,\xi \in \Xi,y \in \ball^n(x,\epsilon_x)$, and $\nu \in \nats$; \label{item:prop:p-env-(a)}
		\item the functions $g^\nu(\cdot, x)$ are continuous on $\Xi$ for any $x$. Indeed, for any $\xi_1,\xi_2 \in \Xi$, we have
		\[
		|g^\nu(\xi_1,x) - g^\nu(\xi_2,x)|\leq L_{\beta,x}^\nu\max\{\|\xi_1\|_2,\|\xi_2\|_2,1\}^{\beta-1}\|\xi_1 - \xi_2\|_2,
		\]
		where $
L_{\beta,x}^\nu=\tfrac{3^{\beta-1}}{\theta^\nu}\max\big\{(2\beta M_x\theta^\nu)^{(\beta-1)/\beta},1\big\}$; \label{item:prop:p-env-(b)}
		\item for any $\xi^\nu\to \xi$ in $\Xi$ and $x^\nu\to x$, $\liminf g^\nu(\xi^\nu, x^\nu) \geq g(\xi, x)$. \label{item:prop:p-env-(c)}
	\end{enumerate}
\end{proposition}

\state Proof.
For Item~\ref{item:prop:p-env-(a)}, let us define functions $l^\nu:\reals^d\times \reals^d \times \reals^n \to \Reals$ as
\[
l^\nu(\zeta,\xi,x)=
\left\{ 
\begin{array}{rl}
g(\zeta,x)+\frac{1}{\beta\theta^\nu}\|\zeta - \xi\|_2^\beta &  \text{if } \zeta,\xi \in \Xi, \\ 
\infty & \text{otherwise}. 
\end{array}\right.
\]
Since $|g(\xi, y)| \leq M_x$ for $\xi \in \Xi,y \in \ball^n(x,\epsilon_x)$ and the set $\Xi$ is closed and nonempty, the functions $(\zeta,\xi,x)\mapsto l^\nu(\zeta,\xi,x)$ are proper, lsc, and level-bounded in $\zeta$ locally uniformly in $(\xi,x)$ in the sense of \cite[Definition 1.16]{VaAn}. Since $g^\nu(\xi,x) = \inf_\zeta l^\nu(\zeta,\xi,x)$ when $\xi \in \Xi$, it follows from \cite[Theorem 1.17]{VaAn} that the functions $g^\nu$ are proper and lsc on $\Xi\times \reals^n$. Also coming from \cite[Theorem 1.17]{VaAn}, for each $(\xi,x) \in \Xi \times \reals^n$, we have $\argmin l^\nu(\cdot,\xi,x) \neq \emptyset$.
 Therefore, for any $x \in \reals^n,y \in \ball^n(x,\epsilon_x)$ and any $\xi \in \Xi$, there exists $\zeta^\nu \in \Xi$ such that 
 \[
 -M_x\leq-M_x+\tfrac{1}{\beta\theta^\nu}\|\xi - \zeta^\nu\|_2^\beta\leq g(\zeta^\nu,y) + \tfrac{1}{\beta\theta^\nu}\|\xi - \zeta^\nu\|_2^\beta = g^\nu(\xi,y)\leq g(\xi,y)\leq M_x.
 \]
 Hence, we get $|g^\nu(\xi, y)| \leq M_x$ for any $x\in\reals^n,\xi \in \Xi,y \in \ball^n(x,\epsilon_x)$, and $\nu \in \nats$.

For Item~\ref{item:prop:p-env-(b)}, let $\xi_1,\xi_2 \in \Xi$. 
There exists $\zeta^\nu \in \Xi$ such that 
\[
-M_x+\tfrac{1}{\beta\theta^\nu}\|\xi_2 - \zeta^\nu\|_2^\beta\leq g(\zeta^\nu,x) + \tfrac{1}{\beta\theta^\nu}\|\xi_2 - \zeta^\nu\|_2^\beta = g^\nu(\xi_2,x)\leq g(\xi_2,x)\leq M_x,
\]
so that $\|\xi_2 - \zeta^\nu \|_2 \leq (2\beta M_x \theta^\nu)^{1/\beta}$. We compute
\begin{align*}
	g^\nu(\xi_1,x) - g^\nu(\xi_2,x)&= g^\nu(\xi_1,x) - g(\zeta^\nu,x) - \tfrac{1}{\beta\theta^\nu}\|\xi_2 - \zeta^\nu\|_2^\beta\\
	&\leq \tfrac{1}{\beta\theta^\nu}\|\xi_1 - \zeta^\nu\|_2^\beta - \tfrac{1}{\beta\theta^\nu}\|\xi_2 - \zeta^\nu\|_2^\beta\\
	&\leq \tfrac{1}{\theta^\nu}\max\{\|\xi_1 - \zeta^\nu\|_2,\|\xi_2 - \zeta^\nu\|_2\}^{\beta-1}\|\xi_1 - \xi_2\|_2 \\
	&\leq \tfrac{1}{\theta^\nu}(\|\xi_1\|_2 + \|\xi_2\|_2 + \|\xi_2 - \zeta^\nu\|_2)^{\beta-1}\|\xi_1 - \xi_2\|_2 \\
	&\leq \tfrac{1}{\theta^\nu}\max\{3\|\xi_1\|_2,3\|\xi_2\|_2,3(2\beta M_x \theta^\nu)^{1/\beta}\}^{\beta-1}\|\xi_1 - \xi_2\|_2,
\end{align*}
where the second inequality uses the triangle inequality
 and the fact that $a^\beta - b^\beta \leq \beta\max\{a,b\}^{\beta-1}|a-b|$ for $a,b\geq 0$, which follows from the mean value theorem applied to $t\mapsto t^\beta$.
Therefore, the claim holds by switching $\xi_1$ and $\xi_2$. 

For Item~\ref{item:prop:p-env-(c)}, let $\xi^\nu \to \xi$ in $\Xi$ and $x^\nu \to x$. For each $\nu \in \nats$, there exists $\zeta^\nu \in \Xi$ such that 
$
 g(\zeta^\nu,x^\nu) + \tfrac{1}{\beta\theta^\nu} \|\zeta^\nu - \xi^\nu\|_2^\beta = g^\nu(\xi^\nu,x^\nu)
$
and $\|\zeta^\nu - \xi^\nu\|_2^\beta\leq 2\beta M_x\theta^\nu$. 
Since $\theta^\nu \to 0$,  we have $\zeta^\nu \to \xi$  and
\[
\nliminf g^\nu(\xi^\nu,x^\nu) =
\liminf g(\zeta^\nu,x^\nu) + \tfrac{1}{\beta\theta^\nu} \|\zeta^\nu - \xi^\nu\|_2^\beta \geq 
\liminf g(\zeta^\nu,x^\nu)  
\geq g(\xi,x),
\]
where the last inequality is from the fact that $g$ is lsc on  $\Xi \times \reals^n$.  
\eop

From Proposition~\ref{prop:p-env}\ref{item:prop:p-env-(b)}, we know that $g^\nu(\cdot, x)$ in \eqref{eq:pth-envelope} are locally Lipschitz continuous on $\Xi$ for any $x$. Therefore, by Propositions~\ref{prop:g-perturb}~and~\ref{prop:p-env}, we obtain the following variational convergence results in the sense of Theorem~\ref{thm:epi-general-impli} for $f^\nu$ in \eqref{eq:general-f-nu}, equipped with the epigraphical regularization in \eqref{eq:pth-envelope}.

\begin{corollary} {\rm (weak convergence).} \label{coro:wk}
Let $g_0: \reals^n \to \Reals$ be a proper, lsc function, and $h: \reals^m \to \Reals$ be a proper, lsc, and nondecreasing function. Assume that the functions $g_i: \Xi \times \reals^n \to \reals$ with $i \in [m]$ are lsc and the functions $g_i(\cdot, x)$ are uniformly bounded on $\Xi$, locally uniformly in $x$. For each $\nu \in \nats$, $i \in [m]$, and $(\xi,x) \in \Xi \times \reals^n$, define
\[
G^\nu(\xi,x)=(g_1^\nu(\xi,x),\ldots,g_m^\nu(\xi, x)),\quad\text{with}\quad
g^\nu_i(\xi, x) = \ninf_{\zeta \in \Xi} g_i(\zeta,x) + \tfrac{1}{\beta\theta^\nu}\|\zeta - \xi\|_2^\beta,
\]
where $\beta \geq 1$. Suppose that $\argmin \phi \neq \emptyset$ with $\phi$ in \eqref{eq:general-problem}. If $\mu^\nu\in \mathscr{P}(\Xi)$ converge to $\mu \in \mathscr{P}(\Xi)$ weakly and $(\theta^\nu,\lambda^\nu) \downarrow 0$ sufficiently slowly, then, for $f$ and $f^\nu$ in \eqref{eq:f} and \eqref{eq:general-f-nu}, Theorem~\ref{thm:epi-general-impli}\ref{item:thm:epi-general-impli-(a)}--\ref{item:thm:epi-general-impli-(f)} hold. 

Specifically, for any $c^\nu\in(0,\infty) \downarrow 0$, there exist $i_\nu^\circ\in \nats \uparrow \infty$ such that for any $i_\nu \in \nats \uparrow \infty$ satisfying $i_\nu \leq i_\nu^\circ$, the sequence $\{(\lambda^\nu,\theta^\nu)\}_\nu$ can be taken as 
\[
\lambda^\nu \in (0,\infty)\to 0,\quad \theta^\nu=c^{i_\nu},\quad \text{and}\quad (\lambda^\nu)^{1/\alpha}i_\nu \to \infty.
\]
A legitimate choice of $c^\nu$ is $c^\nu = \nu^{-1}$.
\end{corollary}
\state Proof. 
Fix $x_0 \in \argmin \phi$. By Proposition~\ref{prop:p-env}\ref{item:prop:p-env-(a)} and \cite[Example 14.31]{VaAn}, the functions $g^\nu_i$ are random lsc and the mappings $G^\nu(\cdot, x)$ are uniformly bounded on $\Xi$, locally uniformly in $x$. We only need to verify that $G^\nu$ and $\lambda^\nu$ satisfy Theorem~\ref{thm:epi-general}\ref{item:thm:epi-general-G}--\ref{item:thm:epi-general-lambda}. 
For Theorem~\ref{thm:epi-general}\ref{item:thm:epi-general-G}, we have $G^\nu(\xi, x) \leq G(\xi, x)$ for any $\xi$ and $x$ by construction. For Theorem~\ref{thm:epi-general}\ref{item:thm:epi-general-liminf}, we aim to use Proposition~\ref{prop:g-perturb}\ref{item:prop:g-perturb-(a)}. To this end, we need to show the functions $g_i^\nu$ with $\theta^\nu=c^\nu$ satisfy Proposition~\ref{prop:g-perturb}\ref{item:lem:g-perturb-cont} and \ref{item:lem:g-perturb-liminf} for each $i \in [m]$, which follows from Proposition~\ref{prop:p-env}\ref{item:prop:p-env-(b)} and \ref{item:prop:p-env-(c)}.
For Theorem~\ref{thm:epi-general}\ref{item:thm:epi-general-lambda}, we apply Proposition~\ref{prop:g-perturb}\ref{item:prop:g-perturb-(b)} and, with a slight abuse of notation, relabel $\theta^\nu=c^{i_\nu}$. This yields $\|\Ex_{\mu^\nu}[  G^\nu(\bm{\xi},x_0) ] - \Ex_\mu[ G^\nu(\bm{\xi},x_0) ]\|_2 \leq \tfrac{\sqrt{m}}{i_\nu}$ for large $\nu$. Hence, if $\lambda^\nu \to 0$ and $(\lambda^\nu)^{1/\alpha}i_\nu \to \infty$, then $(\lambda^\nu)^{-1/\alpha}(\Ex_{\mu^\nu}[  G^\nu(\bm{\xi},x_0) ] - \Ex_\mu[ G^\nu(\bm{\xi},x_0) ]) \to 0$. The claim follows from Theorem~\ref{thm:epi-general-impli}.
\eop

Corollary~\ref{coro:wk} provides an existence results only, as it offers no explicit guidance for selecting $\theta^\nu$ and $\lambda^\nu$. 
After a discussion of setwise convergence,
Section~\ref{sec:metric} presents explicit guidelines derived from analyzing various metrics that quantify the convergence of $\mu^\nu$.

\subsection{Setwise Convergence}\label{sec:setwise}

In the finite distribution setting, Corollary~\ref{coro:finite} shows that the variational convergence in the sense of Theorem~\ref{thm:epi-general-impli} holds by simply setting $G^\nu = G$ for all $\nu \in \nats$. In contrast, for discrete distributions,  Examples~\ref{ex:discreteI} and \ref{ex:discreteII} and Remark~\ref{rmk:alternative} demonstrate that this choice may lead to a failure of convergence when $\mu^\nu$ converge only weakly to $\mu$. 
This subsection generalizes these observations by considering setwise convergence (also known as strong convergence). Recall that $\mu^\nu \in \mathscr{P}(\Xi)$ converge setwise to $\mu \in \mathscr{P}(\Xi)$ if $\mu^\nu(E) \to \mu(E)$ for every $E\in\mathcal{B}(\Xi)$. Analogous to Proposition~\ref{prop:g-perturb}, which addresses weak convergence, we establish a result under setwise convergence.

\begin{proposition}\label{lem:g-perturb-strong}
 Suppose that $\mu^\nu \in \mathscr{P}(\Xi)$ converge to $\mu \in \mathscr{P}(\Xi)$ setwise. 
Let $g: \Xi \times \reals^n \to \reals$ be a random lsc function such that $g(\cdot, x)$ is uniformly bounded on $\Xi$, locally uniformly in $x$. One has:
\begin{enumerate}[label=\textnormal{(\alph*)}]
	\item for any $x^\nu \to x$, $\liminf \Ex_{\mu^\nu}[g(\bm{\xi},x^\nu) ] \geq \Ex_\mu[ g(\bm{\xi},x) ]$; \label{item:lem:g-perturb-strong-(a)}
\item for any $x_0$,
$\Ex_{\mu^\nu}[ g(\bm{\xi}, x_0)  ] \to \Ex_\mu[ g(\bm{\xi}, x_0) ]$. \label{item:lem:g-perturb-strong-(b)}
\end{enumerate}	
\end{proposition}
\state Proof. For Item~\ref{item:lem:g-perturb-strong-(a)}, we aim to apply a version of Fatou's lemma for setwise convergent measures (see, e.g., \cite[Theorem 4.1]{feinberg2014fatou}).
To this end, note that for each $x$, there exist $\epsilon_x > 0$ and $M_x < \infty$ such that $|g(\xi, y)| \leq M_x$ for any $\xi\in\Xi$ and $y \in \ball^n(x,\epsilon_x)$. Choosing $g_n(s)$ (defined in \cite[Theorem 4.1]{feinberg2014fatou}) as the constant function $s\mapsto -M_x$, we verify that the inequalities in \cite[(4.2)]{feinberg2014fatou} hold. Then, from \cite[Theorem 4.1]{feinberg2014fatou} and the fact that $x \mapsto g(\xi, x)$ is lsc for any $\xi$, we have $\liminf \Ex_{\mu^\nu}[g(\bm{\xi},x^\nu) ] \geq \Ex_\mu[ \liminf g(\bm{\xi},x^\nu) ]$ $\geq$ $\Ex_\mu[ g(\bm{\xi},x)]$.
 For Item~\ref{item:lem:g-perturb-strong-(b)}, since $\xi \mapsto g(\xi, x_0)$ is bounded and measurable on $\Xi$, the claim follows from a general convergence result under setwise convergence; see, e.g., \cite[p.~232, Proposition 18]{Royden2ed}.
\eop

Therefore, when $\mu^\nu$ converge setwise to $\mu$, it is not necessary to introduce $G^\nu$ in the construction of $f^\nu$. Variational convergence can still be achieved by selecting a slowly vanishing $\lambda^\nu$.

\begin{corollary} {\rm (setwise convergence).} \label{coro:setwise} Suppose that $\argmin \phi \neq \emptyset$. Let $g_0: \reals^n \to \Reals$ be a proper, lsc function, and $h: \reals^m \to \Reals$ be a proper, lsc, and nondecreasing function. Assume that the component functions $g_i: \Xi \times \reals^n \to \reals$ of $G:\Xi\times \reals^n \to \reals^m$ are random lsc functions and $g_i(\cdot, x)$ are uniformly bounded on $\Xi$, locally uniformly in $x$. Let $G^\nu = G$. If $\mu^\nu\in \mathscr{P}(\Xi)$ converge to $\mu \in \mathscr{P}(\Xi)$ setwise and $\lambda^\nu \in (0,\infty)\to 0$ sufficiently slowly, then, for $f$ and $f^\nu$ in \eqref{eq:f} and \eqref{eq:general-f-nu},  Theorem~\ref{thm:epi-general-impli}\ref{item:thm:epi-general-impli-(a)}--\ref{item:thm:epi-general-impli-(f)} hold. 
\end{corollary}
\state Proof. 
Fix $x_0 \in \argmin \phi$. 
Since $G^\nu = G$, Theorem~\ref{thm:epi-general}\ref{item:thm:epi-general-G} holds trivially. Applying Proposition~\ref{lem:g-perturb-strong}\ref{item:lem:g-perturb-strong-(a)}, we obtain Theorem~\ref{thm:epi-general}\ref{item:thm:epi-general-liminf}.  Moreover, by Proposition~\ref{lem:g-perturb-strong}\ref{item:lem:g-perturb-strong-(b)}, we have $\Ex_{\mu^\nu}[ G(\bm{\xi}, x_0)  ] \to \Ex_\mu[ G(\bm{\xi}, x_0)]$.
Hence, to guarantee Theorem~\ref{thm:epi-general}\ref{item:thm:epi-general-lambda}, it is always possible to choose $\lambda^\nu \in (0,\infty)\to 0$ such that  
\[(\lambda^{\nu})^{-1/\alpha}(\Ex_{\mu^\nu}[ G(\bm{\xi}, x_0)  ] - \Ex_\mu[ G(\bm{\xi}, x_0)])\to 0.
\]
 The claim then follows from Theorem~\ref{thm:epi-general-impli}. 
\eop

\subsection{Convergence in Probability Metrics}\label{sec:metric}

In this subsection, we propose parameter selections guided by probability  metrics.
\subsubsection{Bounded Lipschitz, Fortet-Mourier, and Wasserstein}\label{sec:metric-bl}

We begin with a setting analogous to that of Section~\ref{sec:envelope} for weak convergence. The goal is to derive quantitative variants of Corollary~\ref{coro:wk}, along with suitable choices for $\lambda^\nu$ and $\theta^\nu$. To this end, we consider the bounded Lipschitz metric, which metrizes weak convergence; see, e.g., \cite[Theorem 1.12.4]{vanderVaartWellner.96}.

\paragraph{Bounded Lipschitz Metric.}
Consider the following sets of functions on $\Xi$:
\[
\mathscr{F}_{\rm B}(\Xi)=\left\{ f:\Xi\to [-1,1] \right\},
\quad
\mathscr{F}_{\rm L}(\Xi)=\left\{ f:\Xi \to \reals \Bigm| |f(\xi_1) - f(\xi_2)|\leq \|\xi_1 - \xi_2\|_2,\forall \xi_1,\xi_2 \in \Xi\right\}.
\]
The \emph{bounded Lipschitz metric} between $\mu_1,\mu_2 \in \mathscr{P}(\Xi)$ is defined as 
\[
d_{\rm BL}(\mu_1,\mu_2)=\sup_{f\in \mathscr{F}_{\rm B}(\Xi)\cap \mathscr{F}_{\rm L}(\Xi)}\Big| \int_\Xi f(\xi) (\mu_1 - \mu_2)(\dd \xi) \Big|;
\]
see \cite{romisch2003stability}, \cite{romisch1991distribution}, and \cite[Chapter 1.12]{vanderVaartWellner.96}. Since $d_{\rm BL}$ metrizes weak convergence over $\mathscr{P}(\Xi)$, $\mu^\nu \in \mathscr{P}(\Xi)$ converge weakly to $\mu \in \mathscr{P}(\Xi)$ if and only if $d_{\rm BL}(\mu^\nu, \mu) \to 0$. By using the Pasch-Hausdorff partial envelope in Example~\ref{example:PH},  we have the following quantitative variant of Corollary~\ref{coro:wk}.

\begin{proposition}{\rm (bounded Lipschitz metric).}\label{prop:BL}
	Consider the setting of Corollary~\ref{coro:wk} with $\beta=1$. In particular, for each $\nu \in \nats$  and $(\xi,x) \in \Xi \times \reals^n$, define
\[
G^\nu(\xi,x)=(g_1^\nu(\xi,x),\ldots,g_m^\nu(\xi, x)),\quad\text{with}\quad
g^\nu_i(\xi, x) = \ninf_{\zeta \in \Xi} g_i(\zeta,x) + \tfrac{1}{\theta^\nu}\|\zeta - \xi\|_2, i \in [m].
\]
For $\mu, \mu^\nu \in \mathscr{P}(\Xi)$, suppose that $d_{\rm BL}(\mu^\nu,\mu) \to 0$. Let $\lambda^\nu$ and $\theta^\nu$ be such that 
	\[
	\lambda^\nu \in (0,\infty)\to 0, \quad\theta^\nu \in (0,\infty)\to 0,\quad\text{and}\quad \tfrac{1}{\lambda^\nu}\big( \tfrac{1}{\theta^{\nu}} d_{\rm BL}(\mu^\nu,\mu) \big)^\alpha  \to 0.
	\]	
Then, for $f$ and $f^\nu$ in \eqref{eq:f} and \eqref{eq:general-f-nu}, Theorem~\ref{thm:epi-general-impli}\ref{item:thm:epi-general-impli-(a)}--\ref{item:thm:epi-general-impli-(g)} hold and $f^\nu \eto f$.
\end{proposition}
\state Proof.
Fix any $x_0 \in \dom g_0$. By Proposition~\ref{prop:p-env}\ref{item:prop:p-env-(a)} and \cite[Example 14.31]{VaAn}, the functions $g^\nu_i$ are random lsc and the mappings $G^\nu(\cdot, x)$ are uniformly bounded on $\Xi$,  locally uniformly in $x$.
Hence, there exist $\epsilon_{x_0} > 0$ and $M_{x_0} < \infty$ such that $|g_i^\nu(\xi, y)| \leq M_{x_0}$ for any $\xi \in \Xi, y \in \ball^n(x_0,\epsilon_{x_0}), i \in [m]$, and $\nu \in \nats$.
 We only need to verify that $G^\nu$ and $\lambda^\nu$ satisfy Theorem~\ref{thm:epi-general}\ref{item:thm:epi-general-G}--\ref{item:thm:epi-general-lambda}. 
For Theorem~\ref{thm:epi-general}\ref{item:thm:epi-general-G}, we have $G^\nu(\xi, x) \leq G(\xi, x)$ for any $\xi$ and $x$ by definition. Theorem~\ref{thm:epi-general}\ref{item:thm:epi-general-liminf} follows from Proposition~\ref{prop:g-perturb}\ref{item:prop:g-perturb-(a)}, which is applicable due to $\theta^\nu \in (0,\infty)\to 0$, Proposition~\ref{prop:p-env}\ref{item:prop:p-env-(b)}~and~\ref{item:prop:p-env-(c)}.
For Theorem~\ref{thm:epi-general}\ref{item:thm:epi-general-lambda}, by Proposition~\ref{prop:p-env}\ref{item:prop:p-env-(b)} with $\beta=1$, we know that $\xi\mapsto g^\nu_i(\xi, x_0)$ is $\tfrac{1}{\theta^\nu}$-Lipschitz continuous on $\Xi$. Meanwhile, $\xi\mapsto g_i^\nu(\xi, x_0)$ is $M_{x_0}$-bounded. Hence, 
\[
\left\{\max\{\tfrac{1}{\theta^\nu},M_{x_0}\}^{-1} g^\nu_i(\cdot, x_0) \bigm| \nu \in \nats, i \in [m]\right\}\subseteq \mathscr{F}_{\rm B}(\Xi) \cap \mathscr{F}_{\rm L}(\Xi).
\]
 By definition of $d_{\rm BL}(\mu^\nu,\mu)$ and $\theta^\nu \in (0,\infty) \to 0$, for sufficiently large $\nu$, we have
\[
(\lambda^\nu)^{-1/\alpha}\|\Ex_{\mu^\nu}[  G^\nu(\bm{\xi},x_0) ] - \Ex_\mu[ G^\nu(\bm{\xi},x_0) ]\|_2 \leq \sqrt{m}((\lambda^\nu)^{1/\alpha}\theta^\nu)^{-1} d_{\rm BL}(\mu^\nu,\mu)\to 0.
\]
Therefore, Theorem~\ref{thm:epi-general}\ref{item:thm:epi-general-lambda} holds. Since $x_0$ is arbitrary, the claim follows from Theorem~\ref{thm:epi-general-impli}.
\eop

 Proposition~\ref{prop:BL} provides an explicit guideline for setting $\lambda^\nu$ and $\theta^\nu$ in Corollary~\ref{coro:wk}. For example, a legitimate choice is $\alpha=1, \lambda^\nu = \theta^\nu = d_{\rm BL}(\mu^\nu, \mu)^{1/2 - \epsilon}$ for some $\epsilon \in (0, 1/2)$.
However, Proposition~\ref{prop:BL} applies only to the case $\beta = 1$ in Corollary~\ref{coro:wk}. When using a general parameter $\beta \geq 1$ in the epigraphical regularization of $g$, it is useful to recall the $\beta$th-order Fortet-Mourier metric, which is widely used in the stochastic programming literature; see, e.g., \cite{rachev2002quantitative,schultz2000some,romisch2003stability}.

\paragraph{Fortet-Mourier Metric.} Consider the following set of locally Lipschitz continuous functions on $\Xi$:
\[
\mathscr{F}_\beta(\Xi)=\Big\{ f:\Xi\to \reals \Bigm| |f(\xi_1) - f(\xi_2)|\leq \max\{1,\|\xi_1\|_2,\|\xi_2\|_2\}^{\beta-1} \|\xi_1 - \xi_2\|_2,\forall \xi_1,\xi_2 \in \Xi\Big\},
\]
where $\beta\geq 1$.
The $\beta$th-order \emph{Fortet-Mourier metric} between $\mu_1,\mu_2 \in \mathscr{P}_\beta(\Xi)$ is defined as 
\[
d_{{\rm FM(}\beta{\rm )}}(\mu_1,\mu_2)=\sup_{f\in \mathscr{F}_\beta(\Xi) } \Big| \int_\Xi f(\xi) (\mu_1 - \mu_2)(\dd \xi) \Big|.
\]
For distributions $\mu,\mu^\nu \in \mathscr{P}_\beta(\Xi)$, we have $\dFM{\beta}(\mu^\nu, \mu) \to 0$ if and only of  $\mu^\nu$ converge weakly to $\mu$ and the corresponding sequence of $\beta$th absolute moments converges as well; see \cite[Theorem 6.2.1]{Rachev91}. 

\begin{proposition}{\rm (Fortet-Mourier metric).}\label{prop:FM}
	Consider the setting of Corollary~\ref{coro:wk} with $\beta\geq1$. In particular, for each $\nu \in \nats$, and $(\xi,x) \in \Xi \times \reals^n$, define
\[
G^\nu(\xi,x)=(g_1^\nu(\xi,x),\ldots,g_m^\nu(\xi, x)),\quad\text{with}\quad
g^\nu_i(\xi, x) = \ninf_{\zeta \in \Xi} g_i(\zeta,x) + \tfrac{1}{\beta\theta^\nu}\|\zeta - \xi\|_2^{\beta}, i \in [m].
\]
Suppose that $\dFM{{\beta}}(\mu^\nu,\mu)\to 0$ for $\mu, \mu^\nu \in \mathscr{P}_{\beta}(\Xi)$. Let $\lambda^\nu$ and $\theta^\nu$ be such that 
	\[
	\lambda^\nu \in (0,\infty)\to 0, \quad\theta^\nu \in (0,\infty)\to 0,\quad\text{and}\quad \tfrac{1}{\lambda^\nu}\big(\tfrac{1}{\theta} \dFM{{\beta}}(\mu^\nu,\mu)\big)^\alpha \to 0.
	\]	
Then, for $f$ and $f^\nu$ in \eqref{eq:f} and \eqref{eq:general-f-nu}, Theorem~\ref{thm:epi-general-impli}\ref{item:thm:epi-general-impli-(a)}--\ref{item:thm:epi-general-impli-(g)} hold and $f^\nu \eto f$.
\end{proposition}
\state Proof.
Fix any $x_0 \in \dom g_0$. 
The proof for verifying Theorem~\ref{thm:epi-general}\ref{item:thm:epi-general-G}~and~\ref{item:thm:epi-general-liminf} is similar to the proof of Proposition~\ref{prop:BL}.
For Theorem~\ref{thm:epi-general}\ref{item:thm:epi-general-lambda}, by Proposition~\ref{prop:p-env}\ref{item:prop:p-env-(b)}, we have 
\[
\bigg\{\left(\tfrac{3^{\beta-1}}{\theta^\nu}\max\big\{(2\beta M_{x_0}\theta^\nu)^{(\beta-1)/\beta},1\big\}\right)^{-1} g^\nu_i(\cdot, x_0) \biggm| \nu \in \nats, i \in [m]\bigg\}\subseteq \mathscr{F}_{\beta}(\Xi).
\] 
 By definition of $\dFM{\beta}(\mu^\nu,\mu)$ and $\theta^\nu \in (0,\infty) \to 0$, for sufficiently large $\nu$, we have
\[
(\lambda^\nu)^{-1/\alpha}\|\Ex_{\mu^\nu}[  G^\nu(\bm{\xi},x_0) ] - \Ex_\mu[ G^\nu(\bm{\xi},x_0) ]\|_2 \leq 
3^{\beta-1}\sqrt{m}((\lambda^\nu)^{1/\alpha}\theta^\nu)^{-1} \dFM{\beta}(\mu^\nu,\mu)\to 0.
\]
Therefore, Theorem~\ref{thm:epi-general}\ref{item:thm:epi-general-lambda} holds. Since $x_0$ is arbitrary, the claim follows from Theorem~\ref{thm:epi-general-impli}.
\eop

When $\beta=1$, the $\beta$th-order Fortet-Mourier metric coincides with the Wasserstein $1$-distance by the Kantorovich-Rubinstein dual formula; see \cite{villani2009optimal} for further discussion.
\paragraph{Wasserstein Distance.}
The \emph{Wasserstein $1$-distance} between $\mu_1,\mu_2 \in \mathscr{P}_1(\Xi)$ can be written as

\[
d_{\rm W}(\mu_1,\mu_2)=
\sup_{f\in \mathscr{F}_{\rm L}(\Xi) } \Big| \int_\Xi f(\xi) (\mu_1 - \mu_2)(\dd \xi) \Big|
\]
in the dual form.
It is easy to see that $d_{\rm BL}(\mu_1,\mu_2) \leq d_{\rm W}(\mu_1,\mu_2)= d_{{\rm FM(}1{\rm )}}(\mu_1,\mu_2)$. 
Hence, if $\mu^\nu \in \mathscr{P}_1(\Xi)$ converge to $\mu \in \mathscr{P}_1(\Xi)$ in the Wasserstein $1$-distance, then setting $\lambda^\nu$ and $\theta^\nu$ such that 
\[
	\lambda^\nu \in (0,\infty)\to 0, \quad\theta^\nu \in (0,\infty)\to 0,\quad\text{and}\quad \tfrac{1}{\lambda^\nu} \big(\tfrac{1}{\theta^{\nu}} d_{\rm W}(\mu^\nu,\mu) \big)^\alpha\to 0
	\]
is sufficient to ensure the conclusion of Proposition~\ref{prop:BL}.
 A similar argument applies to the Wasserstein $\beta$-distance (see \cite[Definition 6.1]{villani2009optimal}) by leveraging its monotonicity in  $\beta$; see \cite[Remark 6.6]{villani2009optimal}.

\subsubsection{Minimal Information, TV, and Other  Discrepancies}\label{sec:tv}

Now, we turn our attention to probability metrics that yield quantitative variants of Corollary~\ref{coro:setwise}. We begin with the \emph{minimal information} (mi) probability (pseudo-)metric \cite[Section 2]{rachev2002quantitative}.

\paragraph{Minimal Information Metric.} 
The definition of the minimal information metric depends on the functions appearing in the Corollary~\ref{coro:setwise}. Consider the following set of functions and distributions:
\begin{align*}
\mathscr{F}_{\rm mi}(\Xi)&= \left\{g_i(\cdot, x):\Xi \to \reals \mid x \in \dom g_0, i \in [m]\right\},\\
\mathscr{P}_{\rm mi}(\Xi)&=\bigg\{\mu \in \mathscr{P}(\Xi)\biggm|  \sup_{f \in\mathscr{F}_{\rm mi}(\Xi)} \Big| \int_\Xi f(\xi) \mu(\dd \xi) \Big| < \infty \bigg\},
\end{align*}
where the functions $g_i$ with $i \in \{0,\ldots,m\}$ are defined in the Corollary~\ref{coro:setwise}.
The minimal information metric between $\mu_1,\mu_2 \in \mathscr{P}_{\rm mi}(\Xi)$ is defined as 
\[
d_{\rm mi}(\mu_1, \mu_2)=\sup_{f \in\mathscr{F}_{\rm mi}(\Xi)} \Big| \int_\Xi f(\xi) (\mu_1 - \mu_2)(\dd \xi) \Big|=\sup_{x \in \dom g_0} \max_{i \in [m]} \Big|\int_\Xi g_i(\xi,x)(\mu_1 - \mu_2)(\dd \xi)\Big|.
\]
Notably, for functions $g_i(\cdot, x)$ that are uniformly bounded in $\Xi$, locally uniformly in $x$, the set $\mathscr{P}_{\rm mi}(\Xi)$ coincides with the one defined in \cite[Section 2]{rachev2002quantitative}. This is because, for each $r > 0$, we have
\[
\max_{i \in [m]} \int_\Xi \inf_{x \in \ball^n(0,r) \cap (\dom g_0)} g_i(\xi,x)\mu(\dd \xi) \geq - \max_{y \in \mathcal{G}} M_y > -\infty,
\]
where $\mathcal{G}$ is a finite subset of $\ball^n(0,r) \cap (\dom g_0)$ such that $\ball^n(0,r) \cap (\dom g_0) \subseteq \bigcup_{y \in \mathcal{G}}\ball^n(y,\epsilon_{y})$.

The minimal information metric is widely used in the literature, particularly in works concerning the stability of stochastic programs; see, e.g., \cite{romisch1991distribution,rachev2002quantitative,henrion1999metric,guo2017convergence}. For convergence of distributions in minimal information metric, we have the following quantitative variant of Corollary~\ref{coro:setwise}.

\begin{proposition}{\rm (minimal information metric).}\label{prop:mi}
	Consider the setting of Corollary~\ref{coro:setwise} with $G^\nu = G$ for any $\nu \in \nats$.
	For $\mu, \mu^\nu \in \mathscr{P}_{\rm mi}(\Xi)$, suppose that $d_{\rm mi}(\mu^\nu,\mu)\to 0$. Let $\lambda^\nu$  be such that 
	\[
	\lambda^\nu \in (0,\infty) \to 0, \quad\text{and}\quad \tfrac{1}{\lambda^\nu} d_{\rm mi}(\mu^\nu,\mu)^\alpha \to 0.
	\]	
Then, for $f$ and $f^\nu$ in \eqref{eq:f} and \eqref{eq:general-f-nu}, Theorem~\ref{thm:epi-general-impli}\ref{item:thm:epi-general-impli-(a)}--\ref{item:thm:epi-general-impli-(g)} hold and $f^\nu \eto f$.
\end{proposition}
\state Proof. 
Fix any $x_0 \in \dom g_0$. 
Since $G^\nu = G$, Theorem~\ref{thm:epi-general}\ref{item:thm:epi-general-G} holds trivially. Applying \cite[Proposition 2.1]{rachev2002quantitative}, we obtain Theorem~\ref{thm:epi-general}\ref{item:thm:epi-general-liminf}. For Theorem~\ref{thm:epi-general}\ref{item:thm:epi-general-lambda}, by the definition of $d_{\rm mi}(\mu^\nu,\mu)$, we have  
\[
(\lambda^{\nu})^{-1/\alpha}\|\Ex_{\mu^\nu}[ G(\bm{\xi}, x_0)  ] - \Ex_\mu[ G(\bm{\xi}, x_0)]\|_2 \leq \sqrt{m}(\lambda^{\nu})^{-1/\alpha} d_{\rm mi}(\mu^\nu,\mu) \to 0.
\]
 The claim then follows from Theorem~\ref{thm:epi-general-impli}.
\eop

Setwise convergence on $\mathscr{P}(\Xi)$ is, in general, not metrizable. However, by imposing additional uniformity in the definition of setwise convergence---specifically, requiring uniform convergence over all measurable sets---one recovers the total variation metric.
\paragraph{Total Variation.}  The \emph{total variation} between $\mu_1,\mu_2 \in \mathscr{P}(\Xi)$ is defined as 
\[
d_{\rm TV}(\mu_1, \mu_2)=\sup \Big\{ \int_\Xi f(\xi)(\mu_1 - \mu_2)(\dd \xi) \Bigm|  \text{ measurable }f:\Xi\to [-1,1] \Big\}.
\]
Convergence in total variation leads to another quantitative variant of Corollary~\ref{coro:setwise}.

\begin{proposition}{\rm (total variation).}\label{prop:tv}
	Consider the setting of Corollary~\ref{coro:setwise} with $G^\nu = G$ for any $\nu \in \nats$. For $\mu, \mu^\nu \in \mathscr{P}(\Xi)$, suppose that $d_{\rm TV}(\mu^\nu,\mu)\to 0$. Let $\lambda^\nu$  be such that   
	\[
	\lambda^\nu \in (0,\infty) \to 0, \quad\text{and}\quad \tfrac{1}{\lambda^\nu} d_{\rm TV}(\mu^\nu,\mu)^\alpha \to 0.
	\]	
Then, for $f$ and $f^\nu$ in \eqref{eq:f} and \eqref{eq:general-f-nu}, Theorem~\ref{thm:epi-general-impli}\ref{item:thm:epi-general-impli-(a)}--\ref{item:thm:epi-general-impli-(g)} hold and $f^\nu \eto f$.
\end{proposition}

\state Proof. 
Fix any $x_0 \in \dom g_0$. 
Since $G^\nu = G$, Theorem~\ref{thm:epi-general}\ref{item:thm:epi-general-G} holds trivially. 
Note that the sequence $\{\mu^\nu\}_\nu$ converges to $\mu$  setwise, as $d_{\rm TV}(\mu^\nu,\mu)\to 0$.
Applying Proposition~\ref{lem:g-perturb-strong}\ref{item:lem:g-perturb-strong-(a)}, we obtain Theorem~\ref{thm:epi-general}\ref{item:thm:epi-general-liminf}. For Theorem~\ref{thm:epi-general}\ref{item:thm:epi-general-lambda}, there exists $M_{x_0} \in [0,\infty)$ such that $|g_i(\xi, x_0)| \leq M_{x_0}$ for every $\xi \in \Xi$ and $i \in [m]$.
By the definition of $d_{\rm TV}(\mu^\nu,\mu)$, we have  
\[
(\lambda^{\nu})^{-1/\alpha}\|\Ex_{\mu^\nu}[ G(\bm{\xi}, x_0)  ] - \Ex_\mu[ G(\bm{\xi}, x_0)]\|_2 \leq \sqrt{m}M_{x_0}(\lambda^{\nu})^{-1/\alpha} d_{\rm TV}(\mu^\nu,\mu) \to 0.
\]
 The claim then follows from Theorem~\ref{thm:epi-general-impli}.
\eop

\paragraph{Other  Discrepancies.} 
Total variation is a widely used probability metric in the literature. Several other probability discrepancies---such as the Kullback-Leibler divergence, the Hellinger distance, and the $\mathcal{X}^2$ divergence---are known to upper bound or be closely related to total variation. 
Suppose that $\mu^\nu$ converge to $\mu$ under a probability discrepancy $d(\mu^\nu, \mu)$ that upper bounds $d_{\rm TV}(\mu^\nu, \mu)$. Then, choosing $\lambda^\nu$ satisfying  $\lambda^\nu \in (0,\infty) \to 0$ and $\tfrac{1}{\lambda^\nu} d(\mu^\nu,\mu)^\alpha \to 0$
is sufficient to ensure the conclusion of Proposition~\ref{prop:tv}.
For a detailed comparison of various probability metrics, we refer the reader to \cite{gibbs2002choosing}.

\subsection{Empirical Approximation}\label{sec:empirical}

Among various approximation of a probability distribution $\mu$, the family of empirical measures are of significant importance. Consider a sequence of independent identically distributed (iid) $\Xi$-valued random variables $\{\bm{\xi}^{\nu}\}_\nu$ on a complete probability space $(\Omega,\mathcal{F},\mathbb{P})$ with shared distribution $\mu=\mathbb{P} \circ (\bm{\xi}^\nu)^{-1}$.  For any $\omega \in \Omega$, consider the empirical measures $\mu^\nu(\omega)=\tfrac{1}{\nu}\sum_{k=1}^\nu \delta_{\bm{\xi}^{k}(\omega)}$.  We can write the plug-in function $\phi^\nu$ in \eqref{eq:phi-nu} and approximating Rockafellian $f^\nu$ in \eqref{eq:general-f-nu} explicitly as
\begin{align}
	\phi^\nu(x)(\omega) &=  g_0(x) + h\Big( \tfrac{1}{\nu}\nsum_{k=1}^\nu G(\bm{\xi}^{k}(\omega),x) \Big), \label{eq:empirical-phi} \\ 
	f^\nu(u,x)(\omega) &= g_0(x) + h\Big(u + \tfrac{1}{\nu}\nsum_{k=1}^\nu G^\nu(\bm{\xi}^{k}(\omega),x) \Big) +  \tfrac{1}{\alpha{\lambda^\nu}}  \|u\|_2^\alpha.\label{eq:emprical-f-nu}
\end{align}
  The following example shows that the instability issue observed in Section~\ref{sec:discussion} persists.
\begin{example}{\rm (empirical approximation).}
	\label{ex:empirical-I} 
Consider the following setting: $n=d=1, m=2,  \Xi=[-1,1], \mu$ is the uniform distribution over $\Xi$, $h=\iota_{(-\infty, 0]^2}, g_1(\xi, x)=\max\{1,\min\{\xi,-1\}\},$ and $g_2(\xi, x)=-g_1(\xi, x)$. Then, the actual $\phi$ simplifies to $\phi(x)=x^2$ with $\inf \phi=0$ and $\argmin \phi = \{0\}$.  
	For any $\nu \in \nats$ and $\omega \in \Omega$, when $\mu$ is replaced by the empirical measures $\mu^\nu(\omega)$,  $\phi^\nu$ in \eqref{eq:empirical-phi} can be written as
	\[
	\phi^\nu(x)(\omega)=
	\left\{ \begin{array}{rl}
x^2 &  \text{if } \omega \in S^\nu= \{\omega \in \Omega\mid\tfrac{1}{\nu}\sum_{k=1}^\nu \bm{\xi}^k(\omega)=0\}, \\ 
\infty &  \text{otherwise.}
\end{array}\right.
	\]
	We know that $\mathbb{P}(\cup_\nu S^\nu)=0$. Hence, $\mathbb{P}$-a.s.,
	we have $
	\ninf \phi^\nu = \infty > 0=\ninf \phi$ and $\argmin \phi^\nu = \emptyset$.
	\end{example}
	
Empirical measures $\mu^\nu(\cdot)$ converge weakly to $\mu$ $\mathbb{P}$-a.s., as established by \cite[Theorem 11.4.1]{Dudley_2002}. This allows for the application of results from Section~\ref{sec:metric-bl}, contingent upon an estimation of the convergence rate in a suitable metric; see, e.g., \cite{fournier2015rate}.
Nevertheless, the empirical measures $\mu^\nu(\cdot)$ are more structured compared to general weakly convergent distributions. Leveraging this distinction, we obtain the following result.

\begin{proposition}\label{prop:empirical}
	{\rm (empirical approximation).}
	Suppose that $\argmin \phi \neq \emptyset$. Let $g_0: \reals^n \to \Reals$ be a proper, lsc function, and $h: \reals^m \to \Reals$ be a proper, lsc, and nondecreasing function. For random lsc functions $g_i: \Xi \times \reals^n \to \reals$ with $i \in [m]$, assume that $g_i(\cdot, x)$ are uniformly bounded on $\Xi$, locally uniformly in $x$. Let $\mu^\nu(\cdot)$ be the empirical measures of $\mu$, and suppose $G^\nu = G$. 	
	Let $\lambda^\nu$ be such that
	\[
	\lambda^\nu \in (0,\infty) \to 0, \quad\text{and}
	 \quad   \frac{(\lambda^\nu)^{2/\alpha}\nu}{\log \log (\nu)} \to \infty.
	\]	
Then, for $f$ and $f^\nu$ in \eqref{eq:f} and  \eqref{eq:emprical-f-nu}, $\mathbb{P}$-a.s., Theorem~\ref{thm:epi-general-impli}\ref{item:thm:epi-general-impli-(a)}--\ref{item:thm:epi-general-impli-(g)} hold and $f^\nu \eto f$.
\end{proposition}

\state Proof.
Fix any $x_0 \in \reals^n$. 
Since $G^\nu = G$, Theorem~\ref{thm:epi-general}\ref{item:thm:epi-general-G} holds trivially. 
Note that $\mu^\nu(\cdot)$ converge weakly to $\mu$  $\mathbb{P}$-a.s.~\cite[Theorem 11.4.1]{Dudley_2002}.
By uniform boundedness, using an epigraphical strong law of large number  \cite[Theorem 2.3]{artstein1995consistency},  we obtain Theorem~\ref{thm:epi-general}\ref{item:thm:epi-general-liminf}. For Theorem~\ref{thm:epi-general}\ref{item:thm:epi-general-lambda}, note that $g_i(\cdot, x_0)$ are uniformly bounded on $\Xi$.
By the law of the iterated logarithm (see \cite[Theorem 8.5.2]{durrett2019probability}), we have 
\[
\limsup{} \frac{\sqrt{\nu}\|\Ex_{\mu^\nu(\cdot)}[ G(\bm{\xi}, x_0)  ] - \Ex_\mu[ G(\bm{\xi}, x_0)]\|_2}{\sqrt{\log \log (\nu)}} < \infty, \quad \text{$\mathbb{P}$-almost surely.}
\]
 The claim then follows from Theorem~\ref{thm:epi-general-impli}.
\eop

The proposition shows that the empirical measures $\mu^\nu(\cdot)$ are ``benign'' approximations of $\mu$. Although $\mu^\nu(\cdot)$ may not converge to $\mu$ in minimal information metric, total variation, or setwise, we can still obtain convergence results without leveraging $G^\nu$ due to epigraphical laws of large numbers.
\section{Chance-Constrained Programs}\label{sec:cc}

In this section, we show how our general convergence results can be applied to chance-constrained programs under distributional perturbations. Specifically,
consider the following problem:
\begin{equation}\label{eq:cc-problem}
\nnmin_{x \in \reals^n}\  g_0(x)  \quad \sothat \quad \ \mu(H_i(x)) \geq b_i, \forall i \in [m],
\end{equation}
where $H_i:\reals^n \rightrightarrows \Xi \subseteq \reals^d$ are osc, $b_i \in [0,1]$, and $\mu \in \mathscr{P}(\Xi)$.
A concrete representation of the mappings $H_i$ includes the set defined by inequalities $H_i(x)=\{\xi \in \Xi \mid \ell_i(\xi, x) \leq 0\}$ for some functions $\ell_i: \Xi \times \reals^n \to \reals$. We refer the reader to \cite{romisch1991stability}, \cite{romisch1991distribution}, and \cite{van2020discussion} for more examples and discussions.

 We reformulate the problem \eqref{eq:cc-problem} as minimizing the following extended-valued function:
\begin{equation}\label{eq:cc-phi}
\phi(x)=g_0(x)+\nsum_{i=1}^m\iota_{(-\infty, 0]}\big(b_i- \mu( H_i(x))\big).
\end{equation}
The function $\phi:\reals^n \to \Reals$ in \eqref{eq:cc-phi} can be written in form of the composite formulation \eqref{eq:general-problem} by defining the outer function $h=\iota_{(-\infty, 0]^m}$ and the component functions $g_i(\xi, x)=b_i - \mathbf{1}_{H_i(x)}(\xi)$ of the mapping $G$ for each $i \in [m]$.
It is easy to see that $h$ is proper, lsc, and nondecreasing.
From \cite[Proposition 2.1]{van2020discussion}, we know that the functions $g_i$ are lsc and uniformly bounded on $\Xi \times \reals^n$, hence random lsc by \cite[Example~14.31]{VaAn}. We also define a set-valued mapping $M:\reals^m \rightrightarrows \reals^n$ as
\begin{equation}\label{eq:def-M}
M(y)=\big\{x \in \dom g_0 \bigm| \mu ( H_i(x)) \geq b_i-y_i,  i \in [m]\big\}.
\end{equation}
The set $M(y)$ can be viewed as the feasible set of \eqref{eq:cc-problem} with a right-hand side perturbation by vector $y$; hence, one has $\dom \phi = M(0)$.

Solutions and minimum values of problem~\eqref{eq:cc-problem} can exhibit disproportionate sensitivity to small perturbations in the distribution $\mu$ as illustrated by Examples~\ref{ex:finite-I},~\ref{ex:finite-II},~\ref{ex:discreteI}, and \ref{ex:discreteII}. 
Let $\mu^\nu$ be the perturbed distributions approaching the actual one $\mu$.
We apply our general results from Sections~\ref{sec:general}~and~\ref{sec:prob-metric} to \eqref{eq:cc-problem} and demonstrate how the approximating Rockafellians $f^\nu$ in \eqref{eq:general-f-nu} offer an effective alternative for improving the solution stability of chance-constrained programs.

\subsection{Variational Convergence}\label{sec:cc-vc}

Parallel to the constructions in Section~\ref{sec:prob-metric}, we consider two forms of the approximating Rockafellians $f^\nu$.

\subsubsection{Setting (S1): Absence of Approximating Mappings $G^\nu$}
We first consider the construction where invoking approximating $G^\nu$ to $G$ is unnecessary, like the one used in Sections~\ref{sec:setwise} and \ref{sec:tv}. Specifically, for the chance-constrained problem in \eqref{eq:cc-problem}, consider the 
approximating Rockafellians $f^\nu$ defined as follows
\begin{equation}\label{eq:cc-s1}
	f^\nu(u,x)= g_0(x)+ \nsum_{i=1}^m\iota_{(-\infty, 0]}\big(u_i + b_i -\mu^\nu(H_i(x))\big) + \tfrac{1}{\alpha\lambda^\nu}\|u\|^\alpha_2,
\end{equation}
where $\alpha \geq 1$.
When $\mu^\nu$ converge to $\mu$ in the sense of $d_{\rm mi}(\mu^\nu, \mu) \to 0$ (as defined in Section~\ref{sec:tv}), quantitative stability results for minimizing the plug-in function $\phi^\nu$ in \eqref{eq:phi-nu} are provided in \cite{romisch1991distribution,romisch1991stability}. These quantitative results rely on qualitative stability theory for general parametric problems from \cite{klatte1987note,robinson1987local}. 
A crucial assumption in these works is the inner semicontinuity (isc), or even metric regularity, of the mapping $M$ (or $M^{-1}$); see also \cite{rachev2002quantitative,henrion1999metric}. In contrast, our approach avoids this assumption entirely. 
Specifically, by minimizing the approximating Rockafellian $f^\nu$ in \eqref{eq:cc-s1} instead of the plug-in function $\phi^\nu$ in \eqref{eq:phi-nu}, leveraging our convergence results in Section~\ref{sec:tv}, we circumvent the need for $M$ to be isc. 
This is illustrated in a corollary of Proposition~\ref{prop:mi}.

\begin{corollary}\label{coro:cc-mi} For $\phi$ in \eqref{eq:cc-phi}, suppose that $\argmin \phi \neq \emptyset$. Let $g_0: \reals^n \to \Reals$ be a proper, lsc function, and $H_i:\reals^n \rightrightarrows \Xi$ be osc.  
For $\mu,\mu^\nu \in \mathscr{P}(\Xi)$, suppose that
 $d_{\rm mi}(\mu^\nu,\mu)\to 0$. Let
	 $\lambda^\nu$ be such that 
	\[
	\lambda^\nu \in (0,\infty) \to 0, \quad\text{and}\quad \tfrac{1}{\lambda^\nu} d_{\rm mi}(\mu^\nu,\mu)^\alpha\to 0.
	\]	
	If the sequence $\{(u^\nu, x^\nu)\}_\nu$ is bounded and generated by
	$
	(u^\nu, x^\nu) \in \epsilon^\nu\text{-}\argmin f^\nu
	$
	with $\epsilon^\nu \to 0$,
then, for $f^\nu$ in \eqref{eq:cc-s1}, one has:
\begin{enumerate}[label=\textnormal{(\alph*)}]
	\item $f^\nu(u^\nu, x^\nu) \to \inf \phi$;
	\item $
		\nOutLim{} ( \epsilon^\nu\text{-}\nargmin f^\nu) \subseteq \{0\}\times \nargmin \phi$; in particular, $\nOutLim{} \{x^\nu\}_\nu  \subseteq \argmin \phi$.
\end{enumerate}
\end{corollary}
\state Proof.
Since the function $g_i(\xi, x)=b_i - \mathbf{1}_{H_i(x)}(\xi)$ is uniformly bounded over $\reals^d\times \reals^n$ for any $i \in [m]$, it holds that $\mathscr{P}_{\rm mi}(\Xi) = \mathscr{P}(\Xi)$. 
Note that the functions $f^\nu$ are tight as $\{(u^\nu, x^\nu)\}_\nu$ are bounded and $\epsilon^\nu \to 0$.
The claim follows from Proposition~\ref{prop:mi}.
\eop

A sufficient condition for $\{(u^\nu, x^\nu)\}_\nu$ to be bounded is level-boundedness of $g_0$, which implies that the set $\bigcup_\nu(\epsilon\text{-}\argmin f^\nu)$ is bounded for any $\epsilon \in [0,\infty)$ and bounded $\{\lambda^\nu\}_\nu$.

When minimizing a function with multiple arguments, we can always interchanges the order of minimization and optimize the function partially with respect to some arguments first. In our case, it is natural to consider the functions $x \mapsto \inf_u f^\nu(u,x)$, which can be written explicitly as follows:
\begin{equation}\label{eq:cc-phi-nu-s1}
\phi^{\nu}_f(x)=\ninf_u f^\nu(u,x) =g_0(x) +\tfrac{1}{\alpha\lambda^\nu} \nsum_{i=1}^m \max\big\{0, b_i -\mu^\nu(H_i(x)) \big\}^\alpha.
\end{equation}
For any $\epsilon > 0$ and any $x^\nu \in \epsilon\text{-}\argmin \phi^\nu_f$, there exists $u^\nu$ such that $(u^\nu, x^\nu) \in \epsilon\text{-}\argmin f^\nu$. Hence, we can directly minimize the function $\phi^\nu_f$ without explicitly resorting to $f^\nu$. Notably, minimizing $\phi_f^\nu$ with parameters $\lambda^\nu$ resembles the classic penalty method for constrained optimization. The key distinction, however, is that traditional penalty methods are motivated by computational considerations, whereas the function  $\phi_f^\nu$ arises from the need to stabilize the problem under distributional perturbations.
That is, when $\lambda^\nu$ are chosen in accordance with the conditions in Corollary~\ref{coro:cc-mi}, the resulting penalized formulation exhibits robustness to inaccuracies in the distributions $\mu^\nu$.

Convergence results in Corollary~\ref{coro:cc-mi} are related to the qualitative stability results for stochastic programs in the literature; see, e.g., \cite[Proposition 1]{klatte1987note}, \cite[Theorem 3.2]{romisch1991stability}, \cite[Theorem 1]{henrion1999metric}, and the survey \cite[Theorem 5]{romisch2003stability}. The key difference is that, while existing works deal with the convergence of $\inf \phi^\nu$ and $\epsilon^\nu\text{-}\argmin \phi^\nu$ for the plug-in function $\phi^\nu$ in \eqref{eq:phi-nu}, our result in Corollary~\ref{coro:cc-mi} is applicable to the  penalty-type formulation $\phi_f^\nu$ in \eqref{eq:cc-phi-nu-s1}, which is a by-product of partially minimizing the Rockafellians $f^\nu$ in \eqref{eq:cc-s1}.
The advantage of minimizing $\phi_f^\nu$ rather than $\phi^\nu$, at least from a theoretical perspective, is significant. Existing qualitative stability results for $\phi^\nu$ require an isc-type assumption on the mapping $M$, which is shown to be indispensable in \cite[Example 40]{romisch2003stability}; see also \cite[Proposition 1]{klatte1987note} and the discussion in \cite[p.~497]{romisch2003stability}. In contrast, our Corollary~\ref{coro:cc-mi} for $\phi_f^\nu$ only requires mild assumptions and provides similar convergence results as those for solving $\phi^\nu$. As corollaries to results in Section~\ref{sec:metric}, similar results hold for total variation, Kullback-Leibler divergence, and other discrepancies. We omit the details for brevity.

When $\mu^\nu(\cdot)$ are empirical measures of $\mu$, as a corollary of Proposition~\ref{prop:empirical}, we have the following convergence result, which can be viewed as a variant of the sample-average approximation:
\begin{corollary}\label{coro:cc-empirical}
 For $\phi$ in \eqref{eq:cc-phi}, suppose that $\argmin \phi \neq \emptyset$. Let $g_0: \reals^n \to \Reals$ be a proper, lsc function, and  $H_i:\reals^n \rightrightarrows \Xi$ be osc.  Let $\mu^\nu(\cdot)$ be the empirical measures generated by iid samples $\bm{\xi}^\nu(\cdot)$ on $(\Omega, \mathcal{F}, \mathbb{P})$ with shared distribution $\mu=\mathbb{P}\circ (\bm{\xi}^\nu)^{-1}$. 	
	Let $\lambda^\nu$ be such that
	\[
	\lambda^\nu \in (0,\infty) \to 0, \quad\text{and}
	\quad   \frac{(\lambda^\nu)^{2/\alpha}\nu}{\log \log (\nu)} \to \infty.
	\]		
	If the sequence $\{(u^\nu, x^\nu)\}_\nu$ is bounded and generated by
	$
	(u^\nu, x^\nu) \in \epsilon^\nu\text{-}\argmin f^\nu
	$
	with $\epsilon^\nu \to 0$,
then, for $f^\nu$ in \eqref{eq:cc-s1}, $\mathbb{P}$-a.s., one has:
\begin{enumerate}[label=\textnormal{(\alph*)}]
	\item $f^\nu(u^\nu, x^\nu) \to \inf \phi$;
	\item $
		\nOutLim{} ( \epsilon^\nu\text{-}\nargmin f^\nu) \subseteq \{0\}\times \nargmin \phi$; in particular, $\nOutLim{} \{x^\nu\}_\nu \subseteq \argmin \phi$.
\end{enumerate}
\end{corollary}

\subsubsection{Setting (S2): With Approximating Mappings $G^\nu$}
Next, we consider the case where $\mu^\nu$ converge weakly to $\mu$, a setting that has been scarcely explored at a similar level of generality to ours. As in Sections~\ref{sec:weak} and \ref{sec:metric-bl}, we employ the Pasch-Hausdorff partial envelope (see Example~\ref{example:PH}) to construct the approximations $G^\nu$ to $G$. We emphasize that this is merely one of many possibilities, and other results from Section~\ref{sec:metric-bl} can be applied in a similar manner. Specifically, we use the following approximating Rockafellians $f^\nu$:
\begin{equation}\label{eq:cc-s2}
	f^\nu(u,x)= g_0(x)+ \nsum_{i=1}^m\iota_{(-\infty, 0]}\big(u_i + \Ex_{\mu^\nu}[g_i^\nu(\bm{\xi}, x)]\big) + \tfrac{1}{\alpha\lambda^\nu}\|u\|^\alpha_2,
\end{equation}
where $g_i^\nu:\Xi \times \reals^n$ adopts the Pasch-Hausdorff partial envelope and can be explicitly written as 
	\[
g_i^\nu(\xi,x)=\ninf_{\zeta \in \Xi} b_i - \mathbf{1}_{H_i(x)}(\zeta)+\tfrac{1}{\theta^\nu}\|\zeta-\xi\|_2 
= b_i+\min\{0,\tfrac{1}{\theta^\nu}\dist(\xi, H_i(x)) - 1\}.
\]
As a corollary to Proposition~\ref{prop:BL}, we have the following convergence results for $f^\nu$ in \eqref{eq:cc-s2}.
\begin{corollary}\label{coro:cc-s2} For $\phi$ in \eqref{eq:cc-phi}, suppose that $\argmin \phi \neq \emptyset$. Let $g_0: \reals^n \to \Reals$ be a proper, lsc function, and  $H_i:\reals^n \rightrightarrows \Xi$ be osc. 
For $\mu,\mu^\nu \in \mathscr{P}(\Xi)$,
suppose that $d_{\rm BL}(\mu^\nu,\mu) \to 0$.
	Let $\lambda^\nu$ and $\theta^\nu$ be such that 
	\[
	\lambda^\nu \in(0,\infty) \to 0, \quad\theta^\nu  \in(0,\infty) \to 0,\quad\text{and}\quad \tfrac{1}{\lambda^\nu} \big( \tfrac{1}{\theta^{\nu}} d_{\rm BL}(\mu^\nu,\mu) \big)^\alpha \to 0.
	\]	
	If the sequence $\{(u^\nu, x^\nu)\}_\nu$ is bounded and generated by
	$
	(u^\nu, x^\nu) \in \epsilon^\nu\text{-}\argmin f^\nu
	$
	with $\epsilon^\nu \to 0$,
then, for $f^\nu$ in \eqref{eq:cc-s2}, one has:
\begin{enumerate}[label=\textnormal{(\alph*)}]
	\item $f^\nu(u^\nu, x^\nu) \to \inf \phi$;
	\item $
		\nOutLim{} ( \epsilon^\nu\text{-}\nargmin f^\nu) \subseteq \{0\}\times \nargmin \phi$; in particular, $\nOutLim{} \{x^\nu\}_\nu \subseteq \argmin \phi$.
\end{enumerate}
\end{corollary}

Similar to \eqref{eq:cc-phi-nu-s1}, if we partially minimizing the Rockafellian $(u,x)\mapsto f^\nu(u,x)$ over $u$, we get the following explicit formulation:
\begin{equation}\label{eq:cc-phi-nu-s2}
\phi^{\nu}_f(x)=\ninf_u f^\nu(u,x) =g_0(x) + \tfrac{1}{\alpha\lambda^\nu}\nsum_{i=1}^m \max\Big\{0, b_i+\Ex_{\mu^\nu}\big[\min\{0,\tfrac{1}{\theta^\nu}\dist(\xi, H_i(x)) - 1\}\big] \Big\}^\alpha.
\end{equation}
Although $\phi_f^\nu$ in \eqref{eq:cc-phi-nu-s2} is more intricate than the plug-in function $\phi^\nu$ in \eqref{eq:phi-nu},  
the advantage of minimizing $\phi^\nu_f$, instead of $\phi^\nu$, is substantial. As \cite[p.~499]{romisch2003stability} notes, weak convergence of distributions is a mild condition; the plug-in function $\phi^\nu$ needs stronger assumptions to be stable; see \cite[p.~216]{romisch1991distribution} and \cite[Theorem 6, Examples 7 and 8]{romisch2003stability}. In contrast, by switching to $\phi_f^\nu$ rather than insisting on minimizing $\phi^\nu$, Corollary~\ref{coro:cc-s2} yields a convergence result similar to \cite[Theorem 6]{romisch2003stability} under mild assumptions.

\subsection{Rate of Convergence}

This section establishes the rate of convergence for the results presented in Section~\ref{sec:cc-vc}. Similar to existing quantitative results, we require additional assumptions on the feasible set mapping $M:\reals^m \rightrightarrows \reals^n$, the actual  $\mu$, and the set-valued mappings $H_i$. 
Our first assumption, relevant for both \textbf{(S1)} and \textbf{(S2)}, concerns the metric subregularity $M^{-1}$, or equivalently, the calmness of $M$.
\begin{assumption}
{\rm (metric subregularity).}
\label{assumption:metric}
The inverse mapping $M^{-1}:\reals^n \rightrightarrows \reals^m$ associated with the set-valued mapping $M$ defined in \eqref{eq:def-M} is metrically subregular at every $x \in M(0)$ for point $0$ with modulus $\kappa$. That is, for each $x \in M(0)$, there exists $\tau_x > 0$ such that
	\[
\dist(z, M(0)) \leq \kappa \dist\big(0, M^{-1}(z)\cap \ball^m(0,\tau_x)\big), \quad \forall  z \in \ball^n(x,\tau_x).
\]
\end{assumption}

In the literature, to establish the convergence rate, a widely used assumption is the metric regularity of $M^{-1}$ (see, e.g., \cite[Theorem 1(iii)]{henrion1999metric} and \cite[Theorem 2.3(iii)]{rachev2002quantitative}),
which is a stronger requirement than  metric subregularity used in Assumption~\ref{assumption:metric}. Meanwhile, the work \cite{henrion1999metric} provides an extensive discussion on sufficient conditions for the metric regularity of the mapping $M^{-1}$, hence some of them are also applicable to our Assumption~\ref{assumption:metric}. 

When quantify the convergence for setting \textbf{(S2)}, we need another geometric assumption on $\mu$ and $H_i$, which is related to the notion of Minkowski content in geometric measure theory.

\begin{assumption}\label{assumption:m-content}{\rm (finite upper outer-Minkowski content).}
	There exists a dense subset $D_\rho$ of $\ball^n(0,\rho)$ such that for any $i \in [m]$ and $x \in D_\rho$, the upper outer-Minkowski content of the set $H_i(x) \subseteq \reals^d$ with respect to $\mu$ is finite; i.e.,
	\begin{equation}\label{eq:assumption:m-content}
	\limsup_{\epsilon \downarrow 0} \tfrac{1}{\epsilon}\mu\big((H_i(x)+\ball^d(0,\epsilon))\backslash H_i(x)\big) < \infty,\quad \forall x \in D_\rho, i \in [m].
	\end{equation}
\end{assumption}

Roughly speaking, in Assumption~\ref{assumption:m-content}, we require the sets $H_i(x)$ to have finite perimeter (in some sense) with respect to $\mu$ for every $x$ in a dense subset $D_\rho$ of $\ball^n(0,\rho)$. 
The following result gives sufficient conditions to validate Assumption~\ref{assumption:m-content}.

\begin{proposition}\label{prop:m-content} Suppose the following conditions hold.
\begin{enumerate}[label=\textnormal{(\roman*)}]
	\item The distribution $\mu$ has decomposition $\mu=\mu_{\rm ac} + \mu_{\rm d}$, where $\mu_{\rm ac}\ll \mathcal{L}^d$ with $\mathcal{L}^d$-a.s.~bounded Radon-Nikodym derivative $\frac{\dd \mu}{\dd \mathcal{L}^d}$, and $\mu_{\rm d}$ is a discrete distribution with uniformly discrete support; i.e., $\supp(\mu_{\rm d})=\{\xi_k\mid k \in \nats\}$ with $\ball^d(\xi_j,\epsilon')\cap \{\xi_k\mid k \in \nats \} = \{\xi_j\}$ for any $j \in \nats$ and some $\epsilon' > 0$.\label{item:prop:m-content-i}
	\item There exists a dense set $D_\rho \subseteq \ball^n(0,\rho)$ such that for any $i \in [m]$, $x \in D_\rho$, and the set $K=H_i(x)$, at least one of the following conditions holds. \label{item:prop:m-content-ii}
	\begin{enumerate}[label=\textnormal{(ii-\alph*)}]
	\item 
	 The set $K$ is a convex body; i.e., $K$ is nonempty, compact, and convex.\label{item:prop:m-content-iid}
	\item The set $K$ is compact with Lipschitz boundary; i.e., for any $z \in \bdry K$, there exists a neighborhood $U$ of $z$ such that $U \cap K$ coincides with the rotated epigraph of a Lipschitz function from $\reals^{d-1}$ to $\reals$.\label{item:prop:m-content-iia}
	\item The set $K$ is compact and $t$-rectifiable for some $t < d$; i.e., $K$ can be expressed as the image of a compact subset of $\reals^{t}$ under a Lipschitz mapping from $\reals^{t}$ to $\reals^d$.\label{item:prop:m-content-iib}
	\end{enumerate}	
	\end{enumerate}	
	Then, Assumption~\ref{assumption:m-content} holds.
\end{proposition}
\state Proof.
Fix $x \in D_\rho, i \in [m]$, and let $K_\epsilon=K + \ball^d(0,\epsilon)=H_i(x) + \ball^d(0,\epsilon)$.
From \ref{item:prop:m-content-i}, for sufficiently small $\epsilon > 0$, we have $\{\xi_k\mid k \in \nats\} \cap (K_\epsilon\backslash K) = \emptyset$ and
\[
\mu(K_\epsilon\backslash K)=\mu_{\rm ac}(K_\epsilon\backslash K) + \nsum_{k=1}^\infty p_k\bm{1}_{K_\epsilon\backslash K}(\xi_k) = \mu_{\rm ac}(K_\epsilon\backslash K) \leq M\mathcal{L}^d(K_\epsilon\backslash K),
\]
where $|\frac{\dd \mu}{\dd \mathcal{L}^d}(\xi)| \leq M <\infty$ for $\mathcal{L}^d$-a.s.~$\xi$. 
In what follows, we discuss the conditions in \ref{item:prop:m-content-ii} separately.

If \ref{item:prop:m-content-iid} holds, the conclusion directly follows from the Steiner's formula. Note that a ``convex body'' in \cite{schneider2013convex} need not have interior points \cite[p.~8]{schneider2013convex}. By \cite[(4.8) and Theorem 4.2.1]{schneider2013convex}, we have the expansion $\mathcal{L}^d(K_\epsilon)=\mathcal{L}^d(K) + 2\epsilon V_{d-1}(K)+O(\epsilon^2)$, where $V_{d-1}$ is the $(d-1)$-dimensional intrinsic volume, which is a finite positive measure; see \cite[(4.9)]{schneider2013convex} for a definition.
Therefore, we have
$\limsup_{\epsilon \downarrow 0}\tfrac{1}{\epsilon}\mu(K_\epsilon\backslash K) \leq 2MV_{d-1}(K) < \infty$.

If \ref{item:prop:m-content-iia} holds,  by a result for compact $K$ with nonempty interior \cite[Corollary 1]{ambrosio2008outer}, we obtain
\[
\limsup_{\epsilon \downarrow 0}\tfrac{1}{\epsilon}\mu(K_\epsilon\backslash K)\leq M \limsup_{\epsilon \downarrow 0}\tfrac{1}{\epsilon}\mathcal{L}^d(K_\epsilon\backslash K) = M\mathcal{H}^{d-1}(\bdry K) < \infty,
\]
where $\mathcal{H}^{d-1}$ is the Hausdorff $(d-1)$-dimensional measure in $\reals^d$.

If \ref{item:prop:m-content-iib} holds, then $\mathcal{L}^d(K) = 0$ since $t < d$. Due to $\mu_{\rm ac}\ll \mathcal{L}^d$, we have $\mu_{\rm ac}(K)=0$. This implies 
\[
\limsup_{\epsilon \downarrow 0}\tfrac{1}{\epsilon}\mu(K_\epsilon\backslash K)=\limsup_{\epsilon \downarrow 0}\tfrac{1}{\epsilon}\mu_{\rm ac}(K_\epsilon\backslash K)=\limsup_{\epsilon \downarrow 0}\tfrac{1}{\epsilon}\mu_{\rm ac}(K_\epsilon) \leq M \limsup_{\epsilon \downarrow 0}\tfrac{1}{\epsilon}\mathcal{L}^d(K_\epsilon),
\]
which is proportion to the upper $(d-1)$-dimensional Minkowski content of $K$; see \cite[Definition 1]{ambrosio2008outer} and \cite[p.~273, 3.2.37]{federer1996geometric}.
Also from \ref{item:prop:m-content-iib}, the set $K$ is $t$-rectifiable for $t < d$ in the sense of \cite[p.~251, 3.2.14]{federer1996geometric}. By \cite[p.~275, Theorem 3.2.39]{federer1996geometric}, we know that 
\[
\lim_{\epsilon \downarrow 0}\frac{\mathcal{L}^d(K_\epsilon)}{\mathcal{H}^{d-t}(\ball^{d-t}(0,1))\epsilon^{d-t}} = \mathcal{H}^t(K)<\infty,
\]
where the last inequality is from \cite[p.~103, Theorem 7.5]{mattila1999geometry} and the definition of rectifiability. Therefore, since $t < d$, we have
$\limsup_{\epsilon \downarrow 0}\tfrac{1}{\epsilon}\mu(K_\epsilon\backslash K) < \infty$, which is  actually zero if $t < d-1$.
\eop

Examples of distributions satisfying Proposition~\ref{prop:m-content}\ref{item:prop:m-content-i} include continuous distributions, discrete distributions with uniform discrete support, and their finite convex combination.
The convex body condition in Proposition~\ref{prop:m-content}\ref{item:prop:m-content-iid} is readily verifiable and can be incorporated during the modeling stage.
 For Proposition~\ref{prop:m-content}\ref{item:prop:m-content-iia}~and~\ref{item:prop:m-content-iib}, roughly, we require the relative boundary of the sets $H_i(x)$ to be somehow regular and well-behaved in a Lipschitzian sense.

The final main result gives rate of convergence, with a proof deferred to the end of the paper.
\begin{theorem}{\rm (rate of convergence).}\label{thm:BL-rate}
	 For $\rho \in [0,\infty)$ and $\epsilon \in [0,2\rho]$, suppose that $\ball^n(0,\rho) \cap (\argmin \phi) \neq \emptyset$, $ \ball^{m+n}(0,\rho)  \cap(\argmin f^\nu)\neq \emptyset$, and $\inf f^\nu \in [-\rho, \rho - \epsilon]$. Furthermore, assume that $g_0:\reals^n \to \Reals$ is $L$-Lipschitz continuous on $\dom g_0$ and the mapping $M:\reals^m \rightrightarrows \reals^n$ satisfies Assumption~\ref{assumption:metric} with modulus $\kappa$.  
	For $\mu,\mu^\nu \in \mathscr{P}(\Xi)$, we consider the following two settings.
	\begin{enumerate}[leftmargin=3.2em]	
		\item[{\rm\bf (S1)}] The Rockafellians $f^\nu$ are defined in form of \eqref{eq:cc-s1} with 
		\[
		\lambda^\nu \in (0,\infty)\to 0\quad \text{and}\quad \tfrac{1}{\lambda^\nu} d_{\rm mi}(\mu^\nu, \mu)^\alpha \to 0;
		\]
		\item[{\rm\bf (S2)}] The Rockafellians $f^\nu$ are defined in form of \eqref{eq:cc-s2} with
		\[
		\lambda^\nu \in (0,\infty)\to 0,\quad \theta^\nu\in (0,\infty) \to 0,\quad \text{and}\quad \tfrac{1}{\lambda^\nu}\big(\tfrac{1}{\theta^\nu} d_{\rm BL}(\mu^\nu, \mu)\big)^\alpha\to 0.
		\]
		Assumption~\ref{assumption:m-content} holds.
	\end{enumerate}
	Suppose that the sequence $\{(u^\nu,x^\nu)\}_\nu$ is generated by $
	(u^\nu, x^\nu) \in \epsilon^\nu\text{-}\argmin f^\nu $
	with $\epsilon^\nu \to 0$ and  $\|(u^\nu,x^\nu)\|_2 \leq \rho$.
	Then, for any $\gamma \in (0,1)$, there exists $\bar{\nu}\in\nats$ and $C_1,C_2 \in (0,\infty)$ such that 
	\begin{align*}
	|f^\nu(u^\nu,x^\nu) - \inf \phi| &\leq \eta^\nu + \epsilon^\nu, \tag{minimum value}\\
	\dist(x^\nu, (\epsilon^\nu + 2\eta^\nu)\text{-}\argmin \phi) &\leq \eta^\nu, \tag{solution set}\\
	\mu(H_i(x^\nu)) - b_i &\geq  - \eta^\nu, \forall i \in [m], \tag{constraint violation}\end{align*}
	for any $\nu \geq \bar{\nu}$, where
	\begin{enumerate}[leftmargin=3.2em]	
	\item[{\rm\bf (S1)}] 
	$\displaystyle
	\eta^\nu =  C_1\max\left\{\ (\lambda^\nu)^{1/\alpha},\  d_{\rm mi}(\mu^\nu, \mu),\  \tfrac{1}{\lambda^\nu} d_{\rm mi}(\mu^\nu, \mu)^\alpha \ \right\};
	$
	\item[{\rm\bf (S2)}] 
	$\displaystyle
	\eta^\nu =  C_2 \max\left\{\ \gamma,\ 
	(\lambda^\nu)^{1/\alpha},\  \theta^\nu,\  \tfrac{1}{\theta^\nu}d_{\rm BL}(\mu^\nu, \mu),\  \tfrac{1}{\lambda^\nu} \big( \tfrac{1}{\theta^\nu} d_{\rm BL}(\mu^\nu, \mu)\big)^\alpha \ 
	\right\}.
	$
	\end{enumerate}
	Furthermore, the constant $C_1$ is independent of the choice of $\gamma$.
\end{theorem}

Theorem~\ref{thm:BL-rate} provides convergence rate estimates for the results obtained in settings \textbf{(S1)} (see Corollary~\ref{coro:cc-mi}) and  \textbf{(S2)} (see Corollary~\ref{coro:cc-s2}).
	In particular, for setting \textbf{(S1)}, by choosing $\lambda^\nu=d_{\rm mi}(\mu^\nu, \mu)^{\alpha^2/(\alpha+1)}$ and $\epsilon^\nu=d_{\rm mi}(\mu^\nu, \mu)^{\alpha/(\alpha+1)}$, there exists a constant $C_1'$ such that for large $\nu \in \nats$, we have  
	\begin{equation}\label{eq:cc-rate-s1}
	|f^\nu(u^\nu,x^\nu) - \inf \phi| \leq C_1' \epsilon^\nu, \quad \dist\big(x^\nu, (C_1' \epsilon^\nu)\text{-}\argmin \phi\big) \leq C_1' \epsilon^\nu.
	\end{equation}
For setting \textbf{(S2)}, by choosing $\lambda^\nu=d_{\rm BL}(\mu^\nu, \mu)^{\alpha^2/(2\alpha+2)}$, $\theta^\nu = d_{\rm BL}(\mu^\nu, \mu)^{1/2}$, and $\epsilon^\nu = d_{\rm BL}(\mu^\nu, \mu)^{\alpha/(2\alpha+2)}$, for any $\gamma > 0$, there exists a constant $C_2'$ such that for large $\nu \in \nats$, we have 
	\[
	|f^\nu(u^\nu,x^\nu) - \inf \phi| \leq C_2' (\gamma + \epsilon^\nu), \quad \dist\big(x^\nu, (C_2' \gamma + C_2'\epsilon^\nu)\text{-}\argmin \phi\big) \leq C_2' (\gamma + \epsilon^\nu),
	\]
	which characterizes the rate of convergence to the neighborhood of near-minimizers.
	The presence of a positive constant $\gamma > 0$ is attributed to the lack of uniformity in Assumption~\ref{assumption:m-content}. Indeed, under a strengthened assumption, we can set  $\gamma = 0$ and obtain a rate estimation in a similar form of \eqref{eq:cc-rate-s1}.
	
	\begin{remark}{\rm (uniformity and $\gamma=0$ in {\rm \textbf{(S2)}}).}\label{rmk:uniform}
	Suppose that the condition \eqref{eq:assumption:m-content} in Assumption~\ref{assumption:m-content} holds uniformly over all $x \in D_\rho$; i.e.,
	\begin{equation}\label{eq:rmk:gamma-0}
	\sup_{x \in D_\rho}\limsup_{\epsilon \downarrow 0} \tfrac{1}{\epsilon}\mu\big((H_i(x)+\ball^d(0,\epsilon))\backslash H_i(x)\big) < \infty,\quad \forall i \in [m].
	\end{equation}
	Then, a slight (omitted) modification of the proof of Theorem~\ref{thm:BL-rate} allows us to choose $\gamma = 0$ in the setting {\rm \textbf{(S2)}} of Theorem~\ref{thm:BL-rate}. 
	Hence, under \eqref{eq:rmk:gamma-0}, one has
	\[
	\eta^\nu =  C_2 \max\left\{\ 
	(\lambda^\nu)^{1/\alpha},\  \theta^\nu,\  \tfrac{1}{\theta^\nu}d_{\rm BL}(\mu^\nu, \mu),\  \tfrac{1}{\lambda^\nu} \big( \tfrac{1}{\theta^\nu} d_{\rm mi}(\mu^\nu, \mu)\big)^\alpha\ 
	\right\}
	\]
	in {\rm \textbf{(S2)}} with some $C_2 \in (0,\infty)$.
	In particular, by choosing $\lambda^\nu=d_{\rm BL}(\mu^\nu, \mu)^{\alpha^2/(2\alpha+2)}$, $\theta^\nu = d_{\rm BL}(\mu^\nu, \mu)^{1/2}$, and $\epsilon^\nu = d_{\rm BL}(\mu^\nu, \mu)^{\alpha/(2\alpha+2)}$, there exists a constant $C_2'$ such that for large $\nu \in \nats$, we have 
	\[
	|f^\nu(u^\nu,x^\nu) - \inf \phi| \leq C_2' \epsilon^\nu, \quad \dist\big(x^\nu, (C_2' \epsilon^\nu)\text{-}\argmin \phi\big) \leq C_2' \epsilon^\nu.
	\]
	A sufficient condition for the uniform bound in \eqref{eq:rmk:gamma-0} is that there exists $\bar\rho \in [0,\infty)$ such that for any $x \in D_\rho$, the sets $H_i(x) \subseteq \ball^n(0,\bar\rho)$ are convex bodies. This can be seen from Proposition~\ref{prop:m-content}\ref{item:prop:m-content-iid} and the monotonicity of the $(d-1)$-dimensional intrinsic volume $V_{d-1}(\cdot)$; i.e., $V_{d-1}(K) \leq V_{d-1}(K')$ for  convex bodies $K \subseteq K'$, see \cite[(4.9) and Theorem 4.2.1]{schneider2013convex}.
	\end{remark}
	
	Theorem~\ref{thm:BL-rate} can be used to derive the convergence rate of empirical approximation in two different ways. In \cite[Section 4]{rachev2002quantitative}, results and assumptions are proposed to upper bound the minimal information metric along an empirical process, which, combined with Theorem~\ref{thm:BL-rate}, leads to rate estimation for empirical approximation in Corollary~\ref{coro:cc-empirical}.
	Another approach relies on the rate estimation of the Wasserstein distance for empirical process (see, e.g., \cite{fournier2015rate}), which provides an upper bound on the bounded Lipschitz metric and yields the convergence rate for the Rockafellians $f^\nu$ in \eqref{eq:cc-s2}.

	There are extensive works in the literature on establishing the quantitative stability of minimizing the plug-in function $\phi^\nu$ in \eqref{eq:phi-nu}; see, e.g., the survey \cite{romisch2003stability} and references therein. We now compare our Theorem~\ref{thm:BL-rate} with some of the existing results, beginning with setting \textbf{(S1)}. 
	
	As mentioned in Section~\ref{sec:cc-vc}, a key assumption for analyzing the stability of $\phi^\nu$  is the metric regularity of the mapping $M^{-1}$; see, e.g., \cite[Theorem 1(iii)]{henrion1999metric}, \cite[Theorem 2.3(iii)]{rachev2002quantitative}, and \cite[Theorem 39(ii)]{romisch2003stability}. A similar assumption in our Theorem~\ref{thm:BL-rate} is metric subregularity of $M^{-1}$, which is  weaker than metric regularity. Moreover, to estimate the convergence rate of the solution set, a growth condition is required for the function $g_0$; see \cite[Theorem 1(iv)]{henrion1999metric} and \cite[(9) and Theorem 2.4]{rachev2002quantitative}. This growth condition does not appear in our Theorem~\ref{thm:BL-rate}. 
	These two relaxations of assumptions are achieved by minimizing the approximating Rockafellians $f^\nu$, instead of the plug-in function $\phi^\nu$, and an  epigraphical analysis that allows us to derive an upper bound on  $\dist(x^\nu, \delta^\nu\mbox{-}\argmin \phi)$ with vanishing $\delta^\nu$ rather than on $\dist(x^\nu, \argmin \phi)$; see also \cite[Section 4]{Royset.18}. Notably, while a solution of minimizing $f^\nu$ may not be feasible for the original problem of minimizing $\phi$ due to distributional perturbations, the violation can be quantified in terms of vanishing adjustments to the confidence levels $b_i$.
	Although the assumptions and criteria differ, we still compare the obtained rate to \cite[Theorem 1]{henrion1999metric}. When choosing parameters following the suggestions in \eqref{eq:cc-rate-s1}, our convergence rate of the solution set is of order $d_{\rm mi}(\mu^\nu, \mu)^{\alpha/(\alpha+1)}$, which is a $\tfrac{1}{2}-\tfrac{1}{\alpha+1}$ improvement over \cite[Theorem 1]{henrion1999metric}. For the optimal value, under the same setting, our rate is of order $d_{\rm mi}(\mu^\nu, \mu)^{\alpha/(\alpha+1)}$, which is $\tfrac{1}{\alpha+1}$ slower than that in \cite[Theorem 1]{henrion1999metric}.
	
	Regarding setting \textbf{(S2)}, to our knowledge, the convergence rate under similar assumptions is not available in the literature.  When $\mu^\nu$ converge to $\mu$ weakly but not necessarily in minimal information metric, Theorem~\ref{thm:BL-rate} shows that, under an additional mild geometric assumption, the near-minimizers of $f^\nu$ converge to the neighborhood of $\gamma$-near-minimizers of $f$ at a specified rate, where $\gamma > 0$ can be arbitrarily small and can be set to zero under a stronger assumption (see Remark~\ref{rmk:uniform}).

\vspace{.5em}
\state Proof of Theorem~\ref{thm:BL-rate}.
Fix $\rho >0$. 
Since $\ball^{m+n}(0,\rho) \cap (\argmin f^\nu)  \neq \emptyset$ for any $\nu \in \nats$, the functions $f^\nu$ are tight.
Without loss of generality, we may assume that $\kappa \geq 1, L \geq 1$.	
By Assumption~\ref{assumption:metric} and Lebesgue number lemma,
there exists $\tau > 0$ such that for any $x \in M(0)\cap \ball^n(0,\rho) + \ball^n(0,\tau)$ with $\dist(0,M^{-1}(x)) \leq \tau$, it holds
\[
\dist(x, M(0)) \leq \kappa\dist(0,M^{-1}(x)) = \kappa\dist(\Ex_\mu[G(\bm{\xi}, x)],(-\infty,0]^m).
\]
Let $S_\tau=  M(0) + \ball^n(0,\tau)$. Consider the functions $f^\nu_\tau = f^\nu +\iota_{\reals^m\times S_{\tau}}$. By the convergence results (see Proposition~\ref{prop:mi} for \textbf{(S1)} and Proposition~\ref{prop:BL} for \textbf{(S2)}), we know that for both settings:
\[
		\nOutLim{} ( \epsilon^\nu\text{-}\nargmin f^\nu) \subseteq \nargmin f \subseteq  \{0\}\times M(0).
		\]
From \cite[Theorem 4.10(b)]{VaAn}, for large $\nu$, we have  $\ball^{m+n}(0,\rho)\cap(\epsilon^\nu\text{-}\nargmin f^\nu) \subseteq  \reals^m\times S_{\tau}$. This yields 
\[
\ball^{m+n}(0,\rho)\cap(\epsilon^\nu\text{-}\argmin f_{\tau}^\nu) = \ball^{m+n}(0,\rho)\cap(\epsilon^\nu\text{-}\argmin f^\nu), \quad
\inf f_\tau^\nu = \inf f^\nu,
\]
for large $\nu$.
Therefore, we focus on the truncated functions $f^\nu_\tau$ only.
Using \cite[Theorem 6.56]{primer},  we get 
\[
\begin{aligned}
|\inf f^\nu_\tau - \inf \phi| &\leq \hatsetd_\rho(\epi f^\nu_{\tau}, \epi f), \\
\sup{} \Big\{ \sqrt{\|u\|_2^2+\dist(x, \delta\text{-}\argmin \phi)^2} \Bigm| (u,x) \in \ball^{m+n}(0,\rho)\cap(\epsilon^\nu\text{-}\nargmin f^\nu_\tau) \Big\} &\leq 
\hatsetd_\rho(\epi f^\nu_{\tau}, \epi f),
\end{aligned}
\]
provided $\delta > \epsilon^\nu + 2\hatsetd_\rho(\epi f^\nu_{\tau}, \epi f)$, where $\hatsetd_\rho(\epi f^\nu_{\tau}, \epi f)$ is the \emph{truncated Hausdorff distance} between sets $\epi f^\nu_{\tau}, \epi f \subseteq \reals^{m+n+1}$ using the norm 
\[
\max\{\|(u,x) - (\bar{u},\bar{x})\|_2, |t - \bar{t}|\}\quad \text{with}\quad (u,x,t),(\bar{u},\bar{x},\bar{t}) \in \reals^m\times \reals^n\times \reals;
\]
see \cite[Section 6.J]{primer} for definitions.
To estimate $\hatsetd_\rho(\epi f^\nu_{\tau}, \epi f)$, we use the  Kenmochi condition \cite[Proposition 6.58]{primer}. This yields that $\hatsetd_\rho(\epi f^\nu_{\tau}, \epi f) \leq \eta$ for any $\eta \geq 0$ satisfying
\begin{alignat}{2}
 \inf_{\ball^{m+n}((u,x),\eta)} f^\nu_\tau &\leq \max\{ f(u,x),-\rho\} + \eta, \qquad &&\forall (u,x) \in  \ball^{m+n}(0,\rho) \cap (\lev_{\leq \rho} f), \tag{$T_1$}\label{eq:Kenm-T1}\\
 \inf_{\ball^{m+n}((u,x),\eta)} \mathmakebox[\widthof{$f^\nu_\tau$}][l]{f} &\leq \max\{ f^\nu_\tau(u,x),-\rho\} + \eta, \qquad &&\forall (u,x) \in  \ball^{m+n}(0,\rho) \cap (\lev_{\leq \rho} f^\nu_\tau) \tag{$T_2$}.\label{eq:Kenm-T2}
\end{alignat}
We proceed to bound $\hatsetd_\rho(\epi f^\nu_{\tau}, \epi f)$ for settings {\bf (S1)} and {\bf (S2)} separately. 

\textbf{(S1).} 
For \eqref{eq:Kenm-T1}, if $(u,x) \in  \ball^{m+n}(0,\rho) \cap (\lev_{\leq \rho} f)$, then $u = 0$ and $x \in M(0) \subseteq S_\tau$.
Let $\bar{x} =x$, $\bar{u}=\Ex_{\mu}[G(\bm{\xi},x)] - \Ex_{\mu^\nu}[G(\bm{\xi},x)]$, and
\[
\eta_1=\max\{ d_{\rm mi}(\mu^\nu, \mu), \tfrac{1}{\alpha\lambda^\nu}d_{\rm mi}(\mu^\nu, \mu)^\alpha \}.
\]
Note that $\|u - \bar{u}\|_2 \leq \eta_1$. Then, we have $(\bar{u},\bar{x}) \in \ball^{m+n}((u,x),\eta_1)$ and
\[
\begin{aligned}
\ninf_{\ball^{m+n}((u,x),\eta_1)} f^\nu_\tau &\leq f^\nu_\tau(\bar{u},\bar{x})=g_0(x) + \iota_{(-\infty,0]^m}(\Ex_{\mu}[G(\bm{\xi},x)])+ \tfrac{1}{\alpha\lambda^\nu}\|\Ex_{\mu}[G(\bm{\xi},x)]-\Ex_{\mu^\nu}[G(\bm{\xi},x)]\|_2^\alpha \\
&\leq f(u,x) + \tfrac{1}{\alpha\lambda^\nu}\|\Ex_{\mu}[G(\bm{\xi},x)]-\Ex_{\mu^\nu}[G(\bm{\xi},x)]\|_2^\alpha \leq f(u,x) + \eta_1.
\end{aligned}
\]
For \eqref{eq:Kenm-T2}, if $(u,x) \in \ball^{m+n}(0,\rho) \cap (\lev_{\leq \rho} f_\tau^\nu)$, then $x \in S_\tau$ and $-\rho + \tfrac{1}{\alpha\lambda^\nu}\|u\|^\alpha_2\leq g_0(x)+ \tfrac{1}{\alpha\lambda^\nu}\|u\|^\alpha_2 \leq \rho$, which yields $\|u\|_2 \leq (2\alpha\lambda^\nu \rho)^{1/\alpha}$. Let $\bar{u}=0$ and $\bar{x} \in  M(0)$ such that $\|x - \bar{x}\|_2 =  \dist(x, M(0))$, whose existence is due to $x \in S_\tau \cap \ball^n (0,\rho)$ and $\phi$ is lsc and proper. Let $\eta_2=2L\kappa((2\alpha\lambda^\nu \rho)^{1/\alpha} +  d_{\rm mi}(\mu^\nu, \mu))$. Note that $u+\Ex_{\mu^\nu}[G(\bm{\xi}, x)] \leq 0$. We compute
\begin{align*}
\dist(\Ex_\mu[G(\bm{\xi}, x)],(-\infty,0]^m)&\leq \|u+\Ex_{\mu^\nu}[G(\bm{\xi}, x)] - \Ex_\mu[G(\bm{\xi}, x)]\|_2 \\
&\leq \|u\|_2 +\|\Ex_{\mu^\nu}[G(\bm{\xi}, x)] - \Ex_\mu[G(\bm{\xi}, x)]\|_2 \leq (2\alpha\lambda^\nu \rho)^{1/\alpha} +  d_{\rm mi}(\mu^\nu, \mu).
\end{align*}
Due to $\lambda^\nu \to 0$ and $d_{\rm mi}(\mu^\nu, \mu) \to 0$, for large $\nu$, we have $d(0,M^{-1}(x))=\dist(\Ex_\mu[G(\bm{\xi}, x)],(-\infty,0]^m) \leq \tau$ uniformly in $x$.
 Since $x \in M(0)\cap \ball^n(0,\rho) + \ball^n(0,\tau)$,
by metric subregularity, we have
 \[
 \dist(x, M(0)) \leq \kappa(2\alpha\lambda^\nu \rho)^{1/\alpha} + \kappa d_{\rm mi}(\mu^\nu, \mu).
 \]
Therefore, $(\bar{u},\bar{x}) \in \ball^{m+n}((u,x),\eta_2).$ We compute
\begin{align*}
	\ninf_{\ball^{m+n}((u,x),\eta_2)} f &\leq f(\bar{u},\bar{x})=g_0(\bar{x}) \leq g_0(x) + L\|x - \bar{x}\|_2 + \iota_{(-\infty,0]^m}(u+\Ex_{\mu^\nu}[G(\bm{\xi},x)]) + \tfrac{1}{\alpha\lambda^\nu}\|u\|^\alpha_2\\
	&=f^\nu(u,x) + L\dist(x, M(0))\\
	&\leq f^\nu(u,x) + L\kappa\big((2\alpha\lambda^\nu \rho)^{1/\alpha} +  d_{\rm mi}(\mu^\nu, \mu)\big) \leq f^\nu(u,x) + \eta_2.
\end{align*}
In sum, we have $\hatsetd_\rho(\epi f^\nu_{\tau}, \epi f) \leq \max\{\eta_1,\eta_2\}$, hence there exists $C_1 \in (0,\infty)$ such that
\[
\hatsetd_\rho(\epi f^\nu_{\tau}, \epi f) \leq C_1\max\left\{(\lambda^\nu)^{1/\alpha}, d_{\rm mi}(\mu^\nu, \mu), \tfrac{1}{\lambda^\nu} d_{\rm mi}(\mu^\nu, \mu)^\alpha \right\},
\]	
in the setting \textbf{(S1)} as desired.

\textbf{(S2).}
We skip steps and focus only on the non-trivial parts if the the reasoning is similar to that of \textbf{(S1)}.
For \eqref{eq:Kenm-T1}, if $(u,x) \in \ball^{m+n}(0,\rho)\cap (\lev_{\leq \rho} f)$, then $u = 0$ and $x \in\dom \phi \subseteq S_\tau$.
Let $\bar{x} =x$ and $\bar{u}=\Ex_{\mu}[G^\nu(\bm{\xi},x)] - \Ex_{\mu^\nu}[G^\nu(\bm{\xi},x)]$. 
Let 
\[
\eta_1=\max\{ \tfrac{1}{\theta^{\nu}}d_{\rm BL}(\mu^\nu, \mu), \tfrac{1}{\alpha \lambda^\nu} \big( \tfrac{1}{\theta^\nu}d_{\rm BL}(\mu^\nu, \mu) \big)^\alpha \}.
\]
Since $\|G(\xi, x)\|_\infty \leq 1$ for all $\xi \in \Xi$ and $x \in \reals^n$, by Proposition~\ref{prop:p-env}, we have $\|u - \bar{u}\|_2 \leq \tfrac{1}{\theta^{\nu}}d_{\rm BL}(\mu^\nu, \mu)$. Then, we have $(\bar{u},\bar{x}) \in \ball^{m+n}((u,x),\eta_1)$ and
\[
\begin{aligned}
\ninf_{\ball^{m+n}((u,x),\eta_1)} f^\nu_\tau &\leq f^\nu_\tau(\bar{u},\bar{x})=g_0(x) + \iota_{(-\infty,0]^m}(\Ex_{\mu}[G^\nu(\bm{\xi},x)])+ \tfrac{1}{\alpha\lambda^\nu}\|\Ex_{\mu}[G^\nu(\bm{\xi},x)]-\Ex_{\mu^\nu}[G^\nu(\bm{\xi},x)]\|_2^\alpha \\
&\leq f(u,x) + \tfrac{1}{\alpha\lambda^\nu}\big(\tfrac{1}{\theta^\nu}d_{\rm BL}(\mu^\nu, \mu)\big)^\alpha \leq f(u,x) + \eta_1,
\end{aligned}
\]
where the second inequality uses $G^\nu \leq G$. 
For \eqref{eq:Kenm-T2}, if $(u,x) \in \ball^{m+n}(0,\rho) \cap (\lev_{\leq \rho} f^\nu_\tau)$, then $x \in S_\tau$ and $\|u\|_2 \leq (2\alpha\lambda^\nu \rho)^{1/\alpha}$. 
Let $\mathcal{N}_\gamma \subseteq \ball^{m+n}(0,\rho)\cap (\lev_{\leq \rho} f^\nu_\tau) \cap D_\rho$ be a $\gamma$-net of the compact set $\ball^{m+n}(0,\rho)\cap (\lev_{\leq \rho} f^\nu_\tau)$, where $D_\rho$ is defined in Assumption~\ref{assumption:m-content}. Then, there exists $(\widehat{u}, \widehat{x}) \in \mathcal{N}_\gamma$ such that $\|u - \widehat{u}\|_2 \leq \gamma$ and $\|x - \widehat{x}\|_2 \leq \gamma$.
Let $\bar{u}=0$ and $\bar{x} \in  M(0)$ such that $\|\widehat{x} - \bar{x}\|_2 =  \dist(\widehat{x},M(0))$. 
Note that $\widehat{u}+\Ex_{\mu^\nu}[G^\nu(\bm{\xi}, \widehat{x})] \leq 0$ and $\widehat{x} \in S_\tau$. Hence, by metric subregularity, we have
\begin{align*}
\dist(\widehat{x}, M(0)) &\leq  \kappa\|\widehat{u}+\Ex_{\mu^\nu}[G^\nu(\bm{\xi}, \widehat{x})] - \Ex_\mu[G(\bm{\xi}, \widehat{x})]\|_2 \\
&\leq \kappa\|\widehat{u}\|_2 +\kappa\|\Ex_{\mu^\nu}[G^\nu(\bm{\xi}, \widehat{x})] - \Ex_\mu[G^\nu(\bm{\xi}, \widehat{x})]\|_2 +\kappa\|\Ex_{\mu}[G^\nu(\bm{\xi}, \widehat{x})] - \Ex_\mu[G(\bm{\xi}, \widehat{x})]\|_2.
\end{align*}
The last term can be estimated as follows:
\begin{align*}
\|\Ex_{\mu}[G^\nu(\bm{\xi}, \widehat{x})-G(\bm{\xi}, \widehat{x})]\|_2^2 &= \nsum_{i=1}^m \Big(\int_\Xi \min\{0,\tfrac{1}{\theta^\nu}\dist(\xi, H_i(\widehat{x})) - 1\}+\mathbf{1}_{H_i(\widehat{x})}(\xi)\mu({\rm d}\xi)\Big)^2  \\
&\leq \nsum_{i=1}^m \Big(\int_\Xi \left(\mathbf{1}_{H_i(\widehat{x})}-\mathbf{1}_{H_i(\widehat{x})+\ball^d(0,\theta^\nu)}\right)(\xi)\mu({\rm d}\xi)\Big)^2 \\
&\leq m\max_{1\leq i \leq m} \mu\big((H_i(\widehat{x})+\ball^d(0,\theta^\nu))\backslash H_i(\widehat{x})\big)^2.
\end{align*}
By Assumption~\ref{assumption:m-content} and $\mathcal{N}_\gamma \subseteq D_\rho$, for sufficiently large $\nu \in \nats$,
we have
\[
\max_{1\leq i \leq m} \mu\big((H_i(\widehat{x})+\ball^d(0,\theta^\nu))\backslash H_i(\widehat{x})\big)
\leq \max_{y \in \mathcal{N}_\gamma}  \max_{1\leq i \leq m} \mu\big((H_i(y)+\ball^d(0,\theta^\nu))\backslash H_i(y)\big)
\leq C_{\rm mc}\theta^\nu,
\]
where the constant $C_{\rm mc}\in (0,\infty)$ is defined as
\[
C_{\rm mc}=0.001+\max_{y \in \mathcal{N}_\gamma} \max_{1\leq i \leq m}\limsup_{\epsilon \downarrow 0} \tfrac{1}{\epsilon}\mu\big((H_i(y)+\ball(0,\epsilon))\backslash H_i(y)\big),
\]
which is finite due to Assumption~\ref{assumption:m-content}.
Therefore, we have
\begin{align*}
\dist(\widehat{x}, M(0)) &\leq \kappa(\|u\|_2 + \gamma) + \tfrac{\kappa}{\theta^\nu}d_{\rm BL}(\mu^\nu, \mu) + \kappa\|\Ex_{\mu}[G^\nu(\xi, \widehat{x})-G(\xi, \widehat{x})]\|_2\\
&\leq \kappa((2\alpha\lambda^\nu \rho)^{1/\alpha}+\gamma) + \tfrac{\kappa}{\theta^\nu}d_{\rm BL}(\mu^\nu, \mu) + \kappa \sqrt{m}C_{\rm mc}\theta^\nu.
\end{align*}
Let 
$
\eta_2=2L\kappa(\gamma+(2\alpha\lambda^\nu \rho)^{1/\alpha} + \tfrac{1}{\theta^\nu}d_{\rm BL}(\mu^\nu, \mu) +  \sqrt{m}C_{\rm mc}\theta^\nu).
$
Note that 
\[
\|x - \bar{x}\|_2\leq \gamma+\dist(\widehat{x}, M(0))\leq \gamma+ \kappa((2\alpha\lambda^\nu \rho)^{1/\alpha}+\gamma) + \tfrac{\kappa}{\theta^\nu}d_{\rm BL}(\mu^\nu, \mu) + \kappa \sqrt{m}C_{\rm mc}\theta^\nu.
\]
Hence, $(\bar{u},\bar{x}) \in \ball^{m+n}((u,x),\eta_2).$ We compute
\begin{align*}
	\ninf_{\ball^{m+n}((u,x),\eta_2)} f &\leq f(\bar{u},\bar{x})=g_0(\bar{x})  \\
	&\leq g_0(x) + L(\|x - \widehat{x}\|_2 + \|\widehat{x} - \bar{x}\|_2) + \iota_{(-\infty,0]^m}(u+\Ex_{\mu^\nu}[G(\bm{\xi},x)]) + \tfrac{1}{\alpha\lambda^\nu}\|u\|^\alpha_2\\
	&=f^\nu(u,x) + L\gamma + L\dist(\widehat{x}, M(0))\leq f^\nu(u,x) + \eta_2.
\end{align*}
In sum, we have $\hatsetd_\rho(\epi f^\nu_{\tau}, \epi f) \leq \max\{\eta_1,\eta_2\}$. 
Then, there exists $C_2 \in (0,\infty)$ such that 
\[
\hatsetd_\rho(\epi f^\nu_{\tau}, \epi f) \leq C_2 \max\left\{\gamma,
	(\lambda^\nu)^{1/\alpha}, \theta^\nu, \tfrac{1}{\theta^\nu}d_{\rm BL}(\mu^\nu, \mu), \tfrac{1}{\lambda^\nu} \big( \tfrac{1}{\theta^\nu} d_{\rm BL}(\mu^\nu, \mu)\big)^\alpha
	\right\}
\]
in the setting \textbf{(S2)} as desired.
\eop

\smallskip 

\state Acknowledgment. This work is supported in part by the Office of Naval Research under grants N00014-24-1-2318 and N00014-24-1-2492. 

\bibliographystyle{plain}
\bibliography{refs}

\end{document}

%% file: mac.tex
\usepackage{theorem}
\usepackage{enumerate}
\usepackage{array}
\usepackage[reqno]{amsmath}
\usepackage{amssymb}
\usepackage{mathtools}
\usepackage{latexsym}
\usepackage{makeidx}
\usepackage{fancybox}
\input epsf.sty
\usepackage[dvips]{graphicx,color}
\usepackage{showlabels}
\usepackage{subcaption}
\usepackage[dvipsnames]{xcolor}
\usepackage[normalem]{ulem}

\theoremstyle{change}
\newtheorem{proclaim}{PROCLAIM}[section]
\newtheorem{theorem}[proclaim]{Theorem}

\newtheorem{proposition}[proclaim]{Proposition}
\newtheorem{corollary}[proclaim]{Corollary}
\newtheorem{example}[proclaim]{Example}

\newtheorem{assumption}[proclaim]{Assumption}
\newtheorem{remark}[proclaim]{Remark}

\numberwithin{equation}{section}

\outer\def\proclaim #1. #2\par{\medbreak \noindent{\bf#1.\enspace}{\sl#2}\par
  \ifdim\lastskip<\medskipamount
  \removelastskip\penalty55\medskip\fi}
\def\state #1. { \noindent{\bf#1.\enspace}}
\def\algo #1. { \noindent{\bf#1.\enspace}}

\DeclareMathOperator*{\argmin}{argmin}

\DeclareMathOperator{\bdry}{bdry}

\DeclareMathOperator{\dist}{dist}

\DeclareMathOperator{\dom}{dom}

\DeclareMathOperator{\epi}{epi}

\DeclareMathOperator{\lev}{lev}

\DeclareMathOperator{\supp}{supp}

\newcommand{\comp}{\,{\raise 1pt \hbox{$\scriptstyle\circ$}}\,}

\newcommand{\reals}{\mathbb{R}}
\newcommand{\Reals}{\overline{\mathbb{R}}}
\newcommand{\natnums}{{{\rm l} \kern -.13em {\rm N} }}
\newcommand{\nats}{\mathbb{N}}
\newcommand{\snats}{{I\kern -.29em N}}
\newcommand{\rats}{{Q\kern -.64em \raise 1pt \hbox{$\scriptstyle |$}\;\,}}
\newcommand{\srats}
	{{Q\kern -.56em \raise 1.2pt \hbox{$\scriptscriptstyle /$}\,}}
\newcommand{\ints}{Z\kern -.46em Z}
\newcommand{\ball}{\mathbb{B}}
\newcommand{\pluss}{\hskip1pt \raise1pt\vbox{\hrule width6pt \vskip1pt \hrule
                    width6pt} \kern-4pt{\lower1pt\hbox{\vrule height6pt
		    \kern1pt\vrule height6pt}}\hskip5pt}
\newcommand{\eop}
	{\hfill{$\vcenter{\hrule height1pt \hbox{\vrule width1pt height5pt
   	 \kern5pt \vrule width1pt} \hrule height1pt}$} \medskip}

\newcommand{\setd}{{ d \kern -.15em l}}
\newcommand{\hatsetd}{ d \hat{\kern -.15em l }}

\renewcommand{\epsilon}{\varepsilon}
\renewcommand{\phi}{\varphi}

\newcommand{\tto}{\;{\lower 1pt \hbox{$\rightarrow$}}\kern -12pt
           \hbox{\raise 2.5pt \hbox{$\rightarrow$}}\;}
\newcommand{\overto}[1]{\,{\raise 0pt\hbox{$\rightarrow$}}\kern -9pt
     \hbox{\lower 3pt \hbox{$\scriptscriptstyle#1$}}\hskip6pt}
\newcommand{\underto}[1]{\,{\lower 1pt\hbox{$\rightarrow$}}\kern -9pt
     \hbox{\raise 4pt \hbox{$\,\scriptscriptstyle#1$}}\hskip7pt}
\newcommand{\bigoverto}[1]{{\raise 0pt\hbox{$\,\longrightarrow$}}\kern -16pt
     \hbox{\lower 3pt \hbox{$\scriptscriptstyle#1$}}\hskip4pt}
\newcommand{\bigunderto}[1]{\,{\lower 1pt\hbox{$\longrightarrow$}}\kern -16pt
     \hbox{\raise 4pt \hbox{$\,\scriptscriptstyle#1$}}\hskip6pt}
\newcommand{\bigbigto}[2]{\,{\raise 0pt\hbox{$\,\longrightarrow$}}\kern -16pt
     \hbox{\lower 3pt \hbox{$\scriptscriptstyle#2$}}\kern -10pt
     \hbox{\raise 4pt \hbox{$\,\scriptscriptstyle#1$}}\hskip7pt}
\newcommand{\downto}{{\raise 1pt \hbox{$\scriptscriptstyle \,\searrow\,$}}}
\newcommand{\upto}{{\raise 1pt \hbox{$\scriptscriptstyle \,\nearrow\,$}}}

\newcommand{\notimply}
	{\quad\hbox{$\Longrightarrow \kern -14pt {/}$}\hskip6pt\quad}

\newcommand{\lto}{\,{\lower 1pt\hbox{$\rightarrow$}}\kern -10pt
     \hbox{\raise 4pt \hbox{$\, \scriptstyle l$}}\hskip7pt}
\newcommand{\eto}{\,{\lower 1pt\hbox{$\rightarrow$}}\kern -10pt
     \hbox{\raise 4pt \hbox{$\, \scriptstyle e$}}\hskip7pt}
\newcommand{\hto}{\,{\lower 1pt\hbox{$\rightarrow$}}\kern -11pt
     \hbox{\raise 4pt \hbox{$\, \scriptstyle h$}}\hskip7pt}
\newcommand{\pto}{\,{\lower 1pt\hbox{$\rightarrow$}}\kern -11pt
     \hbox{\raise 4.5pt \hbox{$\, \scriptstyle p$}}\hskip7pt}
\newcommand{\cto}{\,{\lower 1pt\hbox{$\rightarrow$}}\kern -11pt
     \hbox{\raise 4pt \hbox{$\, \scriptstyle c$}}\hskip7pt}
\newcommand{\gto}{\,{\lower 1pt\hbox{$\rightarrow$}}\kern -11pt
     \hbox{\raise 4.5pt \hbox{$\, \scriptstyle g$}}\hskip7pt}
\newcommand{\sto}{\,{\lower 1pt\hbox{$\rightarrow$}}\kern -10pt
     \hbox{\raise 4pt \hbox{$\, \scriptstyle s$}}\hskip7pt}
\newcommand{\awto}{\,{\lower 1pt\hbox{$\rightarrow$}}\kern -15pt
     \hbox{\raise 4pt \hbox{$\, \scriptstyle aw$}}\hskip7pt}
\def\Nto{\,{\raise 1pt\hbox{$\rightarrow$}}\kern -13pt
     \hbox{\lower 3pt \hbox{$\, \scriptstyle N$}}\hskip7pt}
\def\Cto{\,{\raise 1pt\hbox{$\rightarrow$}}\kern -14pt
     \hbox{\lower 3pt \hbox{$\, \scriptstyle C$}}\hskip7pt}
\def\fto{\,{\raise 1pt\hbox{$\rightarrow$}}\kern -14pt
     \hbox{\lower 3pt \hbox{$\, \scriptstyle f$}}\hskip7pt}

\newcommand{\low}[1]{{\lower1pt \hbox{$\scriptstyle #1$}}}
\newcommand{\loww}[1]{{\lower2pt \hbox{$\scriptstyle #1$}}}
\newcommand{\high}[1]{{\raise1pt \hbox{$\scriptstyle #1$}}}

\newcommand{\nsum}{\mathop{\sum}\nolimits}

\newcommand{\nInnLim}{\mathop{\rm LimInn}\nolimits}
\newcommand{\nOutLim}{\mathop{\rm Lim\hspace{-0.01cm}Out}\nolimits}
\newcommand{\nLim}{\mathop{\rm Lim}\nolimits}

\newcommand{\nlim}{\mathop{\rm lim}\nolimits}
\newcommand{\nliminf}{\mathop{\rm liminf}\nolimits}
\renewcommand{\liminf}{\mathop{\rm liminf}}
\renewcommand{\limsup}{\mathop{\rm limsup}}
\newcommand{\nlimsup}{\mathop{\rm limsup}\nolimits}
\newcommand{\ninf}{\mathop{\rm inf}\nolimits}

\newcommand{\nnmin}{\mathop{\rm minimize}}

\newcommand{\nargmin}{\mathop{\rm argmin}\nolimits}

\newcommand{\sothat}{\mathop{\rm subject\; to\;}\nolimits}

\newcommand{\lwdy}[2]{\mathrel{\mathop
        {\raisebox{0.1ex}{\null$#1$}}{\hbox{\kern -1.0em
	{\raisebox{-0.8ex}{$\scriptstyle{\;\to #2}$}}}}}}
\newcommand{\lwwdy}[2]{\mathrel{\mathop
        {\raisebox{0.2ex}{\null$#1$}}{\hbox{\kern -1.0em
	{\raisebox{-1.1ex}{$\scriptstyle{\;\to #2}$}}}}}}
\newcommand{\slwdy}[2]{\scriptsize{{\mathrel{\mathop
        {\raisebox{0.1ex}{\null$#1$}}{\hbox{\kern -1.0em
	{\raisebox{-0.8ex}{$\scriptstyle{\;\to #2}$}}}}}}}}
\newcommand{\slwwdy}[2]{\scriptsize{{\mathrel{\mathop
        {\raisebox{0.2ex}{\null$#1$}}{\hbox{\kern -1.0em
	{\raisebox{-1.1ex}{$\scriptstyle{\;\to #2}$}}}}}}}}

\definecolor{lightgray}{gray}{0.75}
\definecolor{myred}{rgb}{0.55,0,0}
\definecolor{myblue}{rgb}{0,0,0.5} 
\definecolor{mygreen}{rgb}{0,0.5,0} 
\definecolor{purple}{rgb}{0.5,0,0.5} 
\definecolor{turq}{rgb}{0,0.805,0.816} 
\definecolor{maroon}{rgb}{0.51,0,0}
\definecolor{MAROON}{rgb}{0.51,0,0}
\definecolor{redor}{rgb}{0.78,0.078,0.078}
\definecolor{dgreen}{rgb}{0,0.3,0}

\newcommand{\Ex}{\mathbb{E}}

\newcommand{\grill}{{\scriptscriptstyle\#}}
\newcommand{\bcdot}{\,{\raise .2ex \hbox{$\centerdot$}}\,}

